\documentclass[12pt]{article}
\usepackage{amssymb,amsmath,amsthm, amsfonts}
\usepackage{graphicx}
\usepackage{epsfig}
\usepackage{tikz}

\theoremstyle{plain}
\newtheorem{theorem}{Theorem}[section]
\newtheorem{lemma}{Lemma}[section]
\newtheorem{proposition}{Proposition}[section]

\theoremstyle{definition}

\newtheorem{remark}{Remark}[section]

\numberwithin{equation}{section}

\begin{document}
\allowdisplaybreaks
\title{\large\bf Global classical small-data solutions for a  three-dimensional Keller--Segel--Navier--Stokes  system modeling coral fertilization}

\author
{\rm Myowin Htwe \\
\it\small School of Mathematics and Statistics,
 Beijing Institute of Technology,\\
\it\small Beijing, 100081, P.R. China\\[2mm]
  \rm Peter Y.~H.~Pang\\
\it\small Department of Mathematics, National University of Singapore,\\
\it\small 10 Lower Kent Ridge Road, Republic of Singapore 119076\\[2mm]
\rm Yifu Wang\thanks{ \it\small Corresponding author. Email: {\tt wangyifu@bit.edu.cn}}\\
\it\small School of Mathematics and Statistics,
 Beijing Institute of Technology,\\
\it\small Beijing, 100081, P.R. China}
\vspace{-6em}
\date{}
\maketitle

\begin{abstract}
We are concerned with the  Keller--Segel--Navier--Stokes  system 
\begin{equation*}
\left\{
\begin{array}{ll}
\rho_t+u\cdot\nabla\rho=\Delta\rho-\nabla\cdot(\rho \mathcal{S}(x,\rho,c)\nabla c)-\rho m, &\!\!   (x,t)\in \Omega\times (0,T),
\\
m_t+u\cdot\nabla m=\Delta m-\rho m,  &\!\!  (x,t)\in \Omega\times (0,T),
\\
c_t+u\cdot\nabla c=\Delta c-c+m,  & \!\!   (x,t)\in \Omega\times (0,T),
\\
u_t+ (u\cdot \nabla) u=\Delta u-\nabla P+(\rho+m)\nabla\phi,\quad \nabla\cdot u=0,  &\!\!  (x,t)\in \Omega\times (0,T)
\end{array}\right.
\end{equation*}
subject to the boundary condition
 $(\nabla\rho-\rho \mathcal{S}(x,\rho,c)\nabla c)\cdot \nu\!\!=\!\nabla m\cdot \nu=\nabla c\cdot \nu=0,  u=0$ in a  bounded smooth domain $\Omega\subset\mathbb R^3$.  
It is shown that the corresponding problem admits a globally classical solution with exponential decay properties   under the hypothesis that 
$\mathcal{S}\in C^2(\overline\Omega\times [0,\infty)^2)^{3\times 3}$ satisfies
$|\mathcal{S}(x,\rho,c)|\leq C_S $ for  some $C_S>0$,  and the initial data satisfy  certain smallness conditions.

\end{abstract}\vspace{-1em}



\section{Introduction}

Chemotaxis, the biased movement of  individuals in response to 
gradients of certain chemicals,
 has a significant effect on pattern formation in numerous biological contexts (see \cite{Bellomo,Hillen,Painter}).  In particular, the chemotaxis  plays an important role in  
 the reproduction of some  invertebrates such as
  corals, anemones and sea urchins. 
 Indeed, there is experimental evidence  that 
eggs can release a chemical which attracts sperms during the process  of coral fertilization
(\cite{Coll1,Coll2,Miller,Riffell,Spehr}).

The important effect of chemotaxis on the  efficiency of
coral fertilization is investigated by Kiselev and Ryzhik (\cite{Kiselev1,Kiselev2})
via the following chemotaxis system (the densities of  egg and sperm gametes are assumed to be identical)
 \begin{equation}\label{1.1}
n_t+U\cdot\nabla n=\Delta n+\chi\nabla\cdot(n \nabla (\Delta)^{-1}n )-\mu n^q ~~~~ \hbox{in}~~\mathbb{R}^N\times (0,T)
\end{equation}
where $n$ represents the density of egg (sperm) gametes, $U$ is a prescribed solenoidal sea
fluid velocity,  and $\chi>0$ denotes the chemotactic sensitivity constant, $\epsilon n^q$
($q\geq 2$) denotes the fertilization phenomenon. For the Cauchy problem in $\mathbb{R}^2$ with initial datum $n(\cdot,0)=n_0$,
the global-in-time
existence of solutions to \eqref{1.1} ($N = 2, 3$) is proved  under the suitable conditions on initial data.
In addition, they showed  that the total mass $$ \int_{\mathbb{R}^2} n(x,t)dx
\rightarrow n_\infty(\chi,n_0,U)~~\hbox{as} ~t\rightarrow\infty $$
with $n_\infty(\chi,n_0,U)>0$ satisfying $n_\infty(\chi,n_0, U)\rightarrow 0$ as $\chi\rightarrow\infty$ in the
case $q > 2$ of supercritical reaction (\cite{Kiselev2}), whereas in the critical  case $q = 2$, the decay rate of $ \int_{\mathbb{R}^2} n(x,t)dx$ is faster than that of $ 1/\log t$ as $t\rightarrow\infty$, and a weaker
 effect of chemotaxis is observed within finite time intervals (\cite{Kiselev1}). Recently, the total mass behavior of solution to \eqref{1.1} is investigated in \cite{JAEWOOK,Cao,WYF1} when the chemical concentration  is governed by a parabolic equation.  In particular,
 the results of \cite{Cao,WYF1} indicate that unlike in the Cauchy problem, the dynamical behavior of  solution to \eqref{1.1} with $q=2$ in the framework  of bounded domains is essentially independent of  the effect from chemotactic cross-diffusion. More precisely, it is shown in \cite{Cao,WYF1} that  whenever $U$ is a bounded and
sufficiently regular solenoidal vector field,   the component $n$ of any non-trivial classical bounded solution to
\begin{equation}\label{1.2}
\left\{
\begin{array}{ll}
n_t+U\cdot\nabla n=\Delta n-\chi\nabla\cdot (n\nabla c)-\mu n^2, & \quad x\in \Omega, t>0,
\\
c_t+U\cdot\nabla c=\Delta c-c+n,   &\quad x\in \Omega, t>0
\end{array}\right.
\end{equation}
decays to zero in
either of the spaces $L^1(\Omega)$ and $L^\infty(\Omega)$, which can be controlled by
appropriate multiples of $1/(t + 1)$ from above and below, respectively.

Experiments indicate that in certain of chemotaxis motion in a liquid environment, the interaction between cells and the surrounding fluid may substantially affect the behavior thereof (\cite{Kiselev3,Lorz1}).  In the style of \cite{Dombrowski,Tuval}, we hence 
suppose that this interaction occurs not
only through transport but possibly also through a buoyancy-driven feedback of sperm (egg) gametes to the fluid velocity. 
Accordingly, it leads to a refinement of \eqref{1.2} in the framework of  chemotaxis--(Navier--)Stokes system
\begin{equation}\label{1.3}
\left\{
\begin{array}{ll}
n_t+u\cdot\nabla n=\Delta n-\nabla\cdot (n\mathcal{S}(x,n,c)\nabla c)-\mu n^2,
\\
c_t+u\cdot\nabla c=\Delta c-c+n,
\\
u_t+\kappa(u\cdot \nabla)u=\Delta u-\nabla P+n\nabla\phi,
\\
 \nabla\cdot u=0.
\end{array}\right.
\end{equation}
 for  the unknown  density  of sperm (egg) gametes $n$, the signal concentration $c$, the fluid velocity $u$ and the associated pressure $P$
 in the physical domain $\Omega\subset\mathbb R^3$. Here the evolution of velocity
$u$ is governed by the incompressible  (Navier)-Stokes equations, in addition, it is driven  by gametes through  buoyant forces within  a gravitational potential $\phi$,
$\phi\in W^{2,\infty}(\Omega)$ and the chemotactic sensitivity tensor $\mathcal{S}(x,n,c)=(s_{ij}(x,n,c))\in  C^2(\overline\Omega\times [0,\infty)^2)$, $i,j\in \{1,2,3\}$, which reflects that the chemotactic migration may
not necessarily be oriented along the gradient of the chemical signal, but may rather involve rotational flux components (see \cite[sec. 4.2.1]{Othmer} or \cite{Xue} for tensor-valued  sensitivities in the chemotaxis system).

In view of mathematical analysis,  the model \eqref{1.3} compounds the known difficulties in the study of the three-dimensional fluid dynamics with the typical intricacies in the study of chemotactic cross-diffusion reinforced  by signal production. In fact, three-dimensional Navier--Stokes equations are yet lacking
complete existence theory,  
particularly the global solvability in classes of suitably regular functions is yet left as an open problem
except in the cases that the initial data are appropriately small (\cite{Wiegner}).
In addition,
it is observed that when
$\mathcal{S}=\mathcal{S}(x,\rho,c)$ is a tensor, the corresponding chemotaxis--fluid system loses the natural energy structure, which plays a key role in the analysis of the scalar-valued
case (\cite{WinklerTAMS,WinklerCPDE,WinklerJFA,WinklerAIHP}).
Despite these challenges, some comprehensive results on the global-boundedness
and large time behavior of solutions are available in the literature (see \cite{Cao1,Li,Liu,Tao1,Wang1,WinklerJFA,Yu} for example). Indeed, by a continuation argument, authors of \cite{Yu} established the global classical solutions of \eqref{1.3} with $\kappa=1,\mu=0$ decaying to ($\bar{n}_0,\bar{n}_0 ,0 $) exponentially with  $\bar{n}_0=\frac 1{|\Omega|}\int_\Omega n_0(x)dx$ if   $
\|n_0\|_{L^{3}(\Omega)},\|\nabla c_0\|_{L^{3}(\Omega)}$ and $\|u_0\|_{L^{3}(\Omega)}$ are
small enough. In particular,  for the 3D chemotaxis--Stokes variant of \eqref{1.3} with $r n-\mu n^2$ instead of
$\mu n^2$ and 
$\mathcal{S}=\chi$ in the $n-$equation,  the existence of global bounded smooth solutions is proved for  appropriately large $\mu>0$ (\cite{Tao1});
  while the corresponding two-dimensional Navier--Stokes variant thereof possesses  a global bounded classical solution
  for arbitrary
  $\mu>0$ (\cite{TaoZAMP}). In addition,  the
latter two works also provide some results on the asymptotic decay of solutions when
$r=0$, which, in the light of  results of \cite{Cao,WYF1}, indeed seems to decay in time like $\frac 1{t+1}$. 
Furthermore, in the very recent paper \cite{WinklerJFA}, Winkler showed that in the delicate 
three-dimensional
setting,  the  Keller--Segel--Navier--Stokes system considered in \cite{TaoZAMP}
possesses at least one globally generalized solution, and that under an explicit condition on the size of $\mu$ this solution approach a spatially homogeneous equilibrium in their first two components.

From a biological point of view,  it is more realistic to  distinguish between  eggs and sperms, and
 it
thereby becomes possible to take into account that only spermatozoids will be affect by chemotactic attraction, whereas the eggs  are governed by random diffusion, fluid transport and degradation upon contact with sperms
during the coral fertilization process (\cite{Espejo3,Espejo1,Kiselev2}).
In addition, the interaction of the gametes and  the ambient  fluid is not negligible.
 The  gametes are assumed to be transported by the fluid, in turn, the motion
of the latter is driven by gametes through  buoyant forces within  a gravitational potential $\phi$.

As an important step toward the comprehensive understanding of the coral fertilization process, we  shall consider  the large time behavior of the egg--sperm chemotaxis--fluid system.  More precisely, this paper is concerned with the following Keller--Segel--Navier--Stokes system in the spatially three--dimensional setting
\begin{equation}\label{1.4}
\left\{
\begin{array}{ll}
\rho_t+u\cdot\nabla\rho=\Delta\rho-\nabla\cdot(\rho \mathcal{S}(x,\rho,c)\nabla c)-\rho m, &\!\!   (x,t)\in \Omega\times (0,T),
\\
m_t+u\cdot\nabla m=\Delta m-\rho m,  &\!\!  (x,t)\in \Omega\times (0,T),
\\
c_t+u\cdot\nabla c=\Delta c-c+m,  & \!\!   (x,t)\in \Omega\times (0,T),
\\
u_t+ (u\cdot \nabla) u=\Delta u-\nabla P+(\rho+m)\nabla\phi,\quad \nabla\cdot u=0,  &\!\!  (x,t)\in \Omega\times (0,T),
\\
\!\!(\nabla\rho-\rho \mathcal{S}(x,\rho,c)\nabla c)\cdot \nu\!\!=\!\nabla m\cdot \nu=\nabla c\cdot \nu=0,  u=0, & \!\!  \!\!  (x,t)\in \partial\Omega\times (0,T),\\
\!\!\!  \rho(x,0)\!=\!\rho_0(x),m(x,0)\!=\!m_0(x),c(x,0)\!=\!c_0(x),u(x,0)\!=\!u_0(x),&\!\!  x\in\Omega,
\end{array}\right.
\end{equation}
where 
 the sperm $\rho$ chemotactically moves toward the higher concentration of the chemical $c$ released by the egg $m$, while the egg $m$ is merely affected by random diffusion, fluid transport and degradation upon contact with the sperm,   $\mathcal{S}=\mathcal{S}(x,\rho,c)$  satisfies
\begin{equation}\label{1.5}
|\mathcal{S}(x,\rho,c)|\leq C_{\mathcal{S}} \quad\textrm{for some}~ C_{\mathcal{S}}>0,
\end{equation}
 and
\begin{equation}\label{1.6}
\left\{
\begin{array}{ll}
\rho_0\in C^0(\overline{\Omega}),~\rho_0\geq0 ~\hbox{and}~ \rho_0\not\equiv0,
\\
m_0\in C^0(\overline{\Omega}),~m_0\geq0 ~\hbox{and}~ m_0\not\equiv0,
\\
c_0\in W^{1,\infty}(\Omega),~c_0\geq0 ~\hbox{and}~ c_0\not\equiv0,
\\
u_0\in D(A^{\beta}) ~\hbox{for all}~ \beta\in(\frac34,1),
\end{array}\right.
\end{equation}
where $A$ denotes the realization of the Stokes operator in $L^2(\Omega)$.

In the context of these assumptions, our main result  can be stated as follows: 

\begin{theorem} \label{TH.1.1}
Suppose that \eqref{1.5} hold and $\int_{\Omega}\rho_0>\int_{\Omega}m_0$. Let $p_0\in(\frac 32, 3)$, $q_0\in(3,\frac{3p_0}{3-p_0})$.
There exists
$\varepsilon>0$ such that for any initial data $(\rho_0,m_0,c_0,u_0)$ fulfilling
\eqref{1.6} as well as
{\setlength\abovedisplayskip{4pt}
\setlength\belowdisplayskip{4pt}
\begin{align*}
\|\rho_0-\rho_\infty\|_{L^{p_0}(\Omega)}<\varepsilon,\quad  \|m_0\|_{L^{q_0}(\Omega)}<\varepsilon,
\quad\| c_0\|_{L^{\infty}(\Omega)}<\varepsilon, \quad\|u_0\|_{L^{3}(\Omega)}<\varepsilon,
\end{align*}}
\eqref{1.4} admits a global classical solution $(\rho,m,c,u,P)$.
In particular, 
for any $\alpha_1\in(0,\min\{\lambda_1,\rho_\infty\})$, $\alpha_2\in(0,\min\{\alpha_1,\lambda_1',1\})$, 
there exist constants $K_i$, $i=1,2,3,4$,  such that for all $t\geq 1 $
{\setlength\abovedisplayskip{4pt}
\setlength\belowdisplayskip{4pt}\begin{align*}
\|m(\cdot,t)\|_{L^\infty(\Omega)}\leq & K_1e^{-\alpha_1 t},
\\
\|\rho(\cdot,t)-\rho_\infty\|_{L^\infty(\Omega)}\leq & K_2e^{-\alpha_1 t},
\\
\|c(\cdot,t)\|_{W^{1,\infty}(\Omega)}\leq & K_3e^{-\alpha_2t},
\\
\|u(\cdot,t)\|_{L^\infty(\Omega)}\leq & K_4 e^{-\alpha_2 t}.
\end{align*}}
Here $\lambda'_1$ is the first eigenvalue of $A$,  $\lambda_1$ is the first nonzero eigenvalue of $-\Delta$ on $\Omega$  under the Neumann boundary condition, and $\rho_\infty=\frac 1{|\Omega|}(\int_{\Omega}\rho_0-\int_{\Omega}m_0)$.
\end{theorem}

As for the case $\int_{\Omega}\rho_0<\int_{\Omega}m_0$, i.e.,
$m_\infty=\frac 1{|\Omega|}(\int_{\Omega}m_0-\int_{\Omega}\rho_0)>0$, we have
\begin{theorem}
Assume that \eqref{1.5} and $\int_{\Omega}\rho_0<\int_{\Omega}m_0$ hold, and let  $p_0\in(2,3)$, $q_0\in(3,\frac{3p_0}{2(3-p_0)})$. 
Then there exists
$\varepsilon>0$ such that for any initial data $(\rho_0,m_0,c_0,u_0)$ fulfilling
\eqref{1.6} as well as
{\setlength\abovedisplayskip{4pt}
\setlength\belowdisplayskip{4pt}
$$\|\rho_0\|_{L^{p_0}(\Omega)}\leq\varepsilon,\quad  \|m_0-m_\infty\|_{L^{q_0}(\Omega)}\leq\varepsilon,
\quad\|\nabla c_0\|_{L^{3}(\Omega)}\leq\varepsilon, \quad\|u_0\|_{L^{3}(\Omega)}\leq\varepsilon,
$$}
\eqref{1.4} admits a global classical solution $(\rho,m,c,u,P)$.
Furthermore,  for any $\alpha_1\!\in\!(0,\min\{\lambda_1,m_\infty,1\})$, $\alpha_2\!\in\!(0,\min\{\alpha_1,\lambda_1'\})$, there exist constants $K_i>0$, $i=1,2,3,4$, such that
\begin{align*}
&\|m(\cdot,t)-m_\infty\|_{L^\infty(\Omega)}\leq K_1e^{-\alpha_1 t},
\\
&\|\rho(\cdot,t)\|_{L^\infty(\Omega)}\leq K_2e^{-\alpha_1 t},
\\
&\|c(\cdot,t)-m_\infty\|_{W^{1,\infty}(\Omega)}\leq K_3e^{-\alpha_2t},
\\
&\|u(\cdot,t)\|_{L^\infty(\Omega)}\leq K_4 e^{-\alpha_2 t}.
\end{align*}
\end{theorem}
\begin{remark}
In our results, we have excluded the case $\int_{\Omega}\rho_0=\int_{\Omega}m_0$. Indeed,  in the light of results of \cite{Cao,WYF1},  algebraical decay rather than exponential decay of the  solutions  is  expected in this case.
\end{remark}
It is noted that a similar result was proved in \cite{LPW}  for the  three-dimensional Stoke variant of \eqref{1.4}. However, as is well-known,  the  nonlinear
convection $(u\cdot \nabla) u$  in the three-dimensional Navier--Stokes equation  may enforce the  spontaneous emergence of singularities in the sense of blow-up with respect to the norm in
$L^\infty(\Omega)$,  we  thereby  subject  the study of classical solutions of \eqref{1.4} to small initial data by an essentially one-step contradiction argument, unlike  that in the
 two-dimensional case (\cite{Espejo1}). 
Moreover,   in comparison with the chemotaxis--fluid system  considered in \cite{Cao1,Yu}, due to
 {\setlength\abovedisplayskip{4pt}
\setlength\belowdisplayskip{4pt}
 $$\|e^{t\Delta}\omega\|_{L^p(\Omega)}\leq C_1\left(1+t^{-\frac 32(\frac1q-\frac1p)}\right)e^{-\lambda_1t}\|\omega\|_{L^q(\Omega)}
 $$}
 for all $\omega\in L^q(\Omega)$ with $\int_\Omega\omega=0 $, 
  $-\rho m$ in the first equation  of \eqref{1.4} gives rise to some difficulty in mathematical analysis despite its dissipative feature. Indeed,
    the core of this argument is to verify that the interval $(0,T)$ on which  solutions enjoy some exponential decay properties
     can be extended to  $(0,\infty)$,  
  which  accordingly requires
 an appropriate  combination of the mass conservation of $\rho(x,t)-m(x,t)$ with the $L^p-L^q$ estimates for the Neumann heat semigroup.

 The plan of this paper is as follows:  In Section 2, we give
a  local  existence  result and some useful estimates. 
In Section 3,  in  the case of $\mathcal{S}$ vanishing on the boundary,  we give the proof of the main results according to either $\int_{\Omega}\rho_0>\int_{\Omega}m_0$ or $\int_{\Omega}\rho_0<\int_{\Omega}m_0$. In the last  section, on the basis of certain a priori estimates,  the  proof of  our main results for the general $\mathcal{S}$ satisfying \eqref{1.5} is realized via an approximation procedure.  

\section{Preliminaries}  

In this section, we provide some preliminary results that will be used in the subsequent
sections. We begin by recalling the important $L^p-L^q$ estimates for the Neumann heat semigroup  on
bounded domains (\cite{Winkler7}).
\vspace{-1em}
\begin{lemma}\label{Lemma 2.1} (Lemma 1.3 of \cite{Winkler7})
Let $(e^{t\Delta})_{t>0}$ denote the Neumann heat semigroup in the domain $\Omega$  and $\lambda_1>0$ denote the first nonzero eigenvalue of $-\Delta$ in $\Omega\subset\mathbb R^N$  under the Neumann boundary condition. There exists $C_i$, $i=1,2,3,4$, such that for all $t>0$,

$(i)$ If $1\leq q\leq p\leq\infty$, then for all $\omega\in L^q(\Omega)$ with $\int_\Omega\omega=0$,
$$\|e^{t\Delta}\omega\|_{L^p(\Omega)}\leq C_1\left(1+t^{-\frac N2(\frac1q-\frac1p)}\right)e^{-\lambda_1t}\|\omega\|_{L^q(\Omega)};$$

$(ii)$ If $1\leq q\leq p\leq\infty$, then for all $\omega\in L^q(\Omega)$,
$$\|\nabla e^{t\Delta}\omega\|_{L^p(\Omega)}\leq C_2\left(1+t^{-\frac12-\frac N2(\frac1q-\frac1p)}\right)e^{-\lambda_1t}\|\omega\|_{L^q(\Omega)};$$

$(iii)$ If $2\leq q\leq p<\infty$, then for all $\omega\in W^{1,q}(\Omega)$,
$$\|\nabla e^{t\Delta}\omega\|_{L^p(\Omega)}\leq C_3\left(1+t^{-\frac N2(\frac1q-\frac1p)}\right)e^{-\lambda_1t}\|\nabla\omega\|_{L^q(\Omega)};$$

$(iv)$ If $1\leq q\leq p\leq\infty$ or $1<q<\infty$ and $p=\infty$, then for all $\omega\in (L^q(\Omega))^N$,
$$\|e^{t\Delta}\nabla\cdot\omega\|_{L^p(\Omega)}\leq C_4\left(1+t^{-\frac12-\frac N2(\frac1q-\frac1p)}\right)e^{-\lambda_1t}\|\omega\|_{L^q(\Omega)}.$$
\end{lemma}

Next we introduce the Stokes operator and recall estimates for the corresponding semigroup.
With $L_\sigma^p(\Omega):=\{\varphi\in L^p(\Omega)|\nabla\cdot\varphi=0\}$ and $\mathcal{P}$ representing the Helmholtz projection of $L^p(\Omega)$ onto $L_\sigma^p(\Omega)$, the Stokes operator on $L_\sigma^p(\Omega)$ is defined as $A_p=-\mathcal{P}\Delta$ with domain
$D(A_p):=W^{2,p}(\Omega)\cap W^{2,p}_0(\Omega)\cap L_\sigma^p(\Omega)$. Since $A_{p_1}$ and $A_{p_2}$ coincide on the intersection of their domains for $p_1$, $p_2\in(1,\infty)$, we will drop the index in the following.

\begin{lemma}\label{Lemma 2.2}(Lemma 2.3 of \cite{Cao1})
The Stokes operator $A$ generates the analytic semigroup $(e^{-tA})_{t>0}$ in $L_\sigma^r(\Omega)$. Its spectrum satisfies $\lambda_1'=\hbox{inf}~ \hbox{Re}\sigma(A)>0$ and we fix $\mu\in(0,\lambda_1')$. For any such $\mu$, we have

$(i)$ For any $p\in(1,\infty)$ and $\gamma\geq0$,  there is $C_5(p,\gamma)>0$ such that for all $\phi\in L^p_\sigma(\Omega)$,
$$\|A^\gamma e^{-tA}\phi\|_{L^p(\Omega)}\leq C_5(p,\gamma) t^{-\gamma}e^{-\mu t}\|\phi\|_{L^p(\Omega)};$$

$(ii)$ For any $p$, $q$ with $1<p\leq q<\infty$, there is $C_6(p,q)>0$ such that for all $\phi\in L^p_\sigma(\Omega)$,
$$\|e^{-tA}\phi\|_{L^q(\Omega)}\leq C_6(p,q) t^{-\frac N2\left(\frac1p-\frac1q\right)}e^{-\mu t}\|\phi\|_{L^p(\Omega)};$$

$(iii)$ For any $p$, $q$ with $1<p\leq q<\infty$,  there is $C_7(p,q)>0$ such that for all $\phi\in L^p_\sigma(\Omega)$,
$$\|\nabla e^{-tA}\phi\|_{L^q(\Omega)}\leq C_7(p,q) t^{-\frac12-\frac N2\left(\frac1p-\frac1q\right)}e^{-\mu t}\|\phi\|_{L^p(\Omega)};$$

$(iv)$ If $\gamma\geq0$ and $1<p<q<\infty$ satisfy $2\gamma-\frac Nq\geq1-\frac Np$, there is $C_8(\gamma,p,q)>0$ such that for all $\phi\in D(A_q^\gamma)$,
$$\|\phi\|_{W^{1,p}(\Omega)}\leq C_8(\gamma,p,q)\|A^\gamma\phi\|_{L^q(\Omega)}.$$
\end{lemma}

\begin{lemma}\label{Lemma2.3} (Theorem 1 and Theorem 2 of \cite{Fujiwara})
The Helmholtz projection $\mathcal{P}$ defines a bounded linear operator $\mathcal{P}$: $L^p(\Omega)\rightarrow L^p_\sigma(\Omega)$; in particular, for any $p\in(1,\infty)$, there exists $C_9(p)>0$ such that
$\|\mathcal{P}\omega\|_{L^p(\Omega)}\leq C_9(p)\|\omega\|_{L^p(\Omega)}$
for every $\omega\in L^p(\omega)$.
 \end{lemma}
The following elementary lemma provides some useful information on both the short-time and
the large-time behavior of certain integrals,  which is used
in the proof of the main results.

\begin{lemma}\label{Lemma 2.4}(Lemma 1.2 of \cite{Winkler7}) 
Let $\alpha\in (0,1),\beta\in (0,1)$, and  $\gamma$, $\delta$  be positive constants such that $ \gamma \neq\delta$. Then there exists $C_{10}(\alpha,\beta,\gamma,\delta) > 0$ such
that
$$\int_0^t(1+s^{-\alpha})(1+(t-s)^{-\beta})e^{-\gamma s}e^{-\delta(t-s)}ds\leq C_{10}(\alpha,\beta,\gamma,\delta)  \left(1+t^{\min\{0,1-\alpha-\beta\}}\right)e^{-\min\{\gamma,\delta\}t}.$$
\end{lemma}

Next we recall the result on the local  existence of classical solutions, which
can be proved by a straightforward adaptation of well-known fixed point argument  
(see \cite{WinklerCPDE} for example).\vspace{-1em}

\begin{lemma}\label{lemma2.5}
Suppose that \eqref{1.5}, \eqref{1.6} and
{\setlength\abovedisplayskip{4pt}
\setlength\belowdisplayskip{4pt}
\begin{equation}\label{2.1}
\mathcal{S}(x,\rho,c)=0, ~~(x,\rho,c)\in \partial\Omega\times [0,\infty)\times [0,\infty)
\end{equation}}
hold. Then there exist $T_{max}\in(0,\infty]$  and a classical solution
$(\rho,m,c,u,P)$ of \eqref{1.4} on $(0,T_{max})$. Moreover, $\rho,m,c$ are nonnegative in $\Omega\times(0,T_{max})$, and
if $T_{max}<\infty$, then for $\beta\in(\frac34,1)$,
$$\|\rho(\cdot,t)\|_{L^\infty(\Omega)}+\|m(\cdot,t)\|_{L^\infty(\Omega)}
+\|c(\cdot,t)\|_{W^{1,\infty}(\Omega)}+
\|A^{\beta}u(\cdot,t)\|_{L^2(\Omega)}\rightarrow \infty ~~\hbox{as}~~t\rightarrow T_{max}.
$$
This solution is unique, up to addition of constants to $P$.\vspace{-0.5em}
\end{lemma}

The following  elementary properties of the solutions in Lemma \ref{lemma2.5} are immediate consequences of the integration of
the first and second equations in \eqref{1.4}, as well as an application of the maximum principle to the second and third equations.\vspace{-1em}

\begin{lemma}\label{Lemma2.6}
Suppose that \eqref{1.5}, \eqref{1.6} and \eqref{2.1} hold. Then for all $t\in(0,T_{max})$, the solution of \eqref{1.4} from Lemma \ref{lemma2.5} satisfies
{\setlength\abovedisplayskip{4pt}
\setlength\belowdisplayskip{4pt}\begin{align}
&\|\rho(\cdot,t)\|_{L^1(\Omega)}\leq\|\rho_0\|_{L^1(\Omega)},\quad\|m(\cdot,t)\|_{L^1(\Omega)}\leq\|m_0\|_{L^1(\Omega)},
\label{2.2}
\\
&\int_0^t\|\rho(\cdot,s)m(\cdot,s)\|_{L^1(\Omega)}ds\leq\min\{\|\rho_0\|_{L^1(\Omega)},\|m_0\|_{L^1(\Omega)}\},
\label{2.3}
\\
&\|\rho(\cdot,t)\|_{L^1(\Omega)}-\|m(\cdot,t)\|_{L^1(\Omega)}=\|\rho_0\|_{L^1(\Omega)}-\|m_0\|_{L^1(\Omega)},
\label{2.4}
\\
&\|m(\cdot,t)\|_{L^2(\Omega)}^2+2\int_0^t\|\nabla m(\cdot,s)\|_{L^2(\Omega)}^2ds\leq \|m_0\|_{L^2(\Omega)}^2,
\label{2.5}
\\
&\|m(\cdot,t)\|_{L^\infty(\Omega)}\leq\|m_0\|_{L^\infty(\Omega)},
\label{2.6}
\\
&\|c(\cdot,t)\|_{L^\infty(\Omega)}\leq\max\{\|m_0\|_{L^\infty(\Omega)},\|c_0\|_{L^\infty(\Omega)}\}.
\label{2.7}
\end{align}}
\end{lemma}
\section {Proof of Theorems  for $\mathcal{S}=0$ on $\partial\Omega$} 

In this section, we give the proofs of Theorem 1.1 and Theorem 1.2 when $\mathcal{S}=0$ on $\partial\Omega$, respectively, i.e. the proof of  Proposition 3.1 and Proposition 3.2 below,
under which the boundary condition for $\rho$ in \eqref{1.4} actually  reduces to the homogeneous
Neumann condition $\nabla \rho\cdot \nu=0$.

In the case $\int_{\Omega}\rho_0>\int_{\Omega}m_0$, i.e.,
$\rho_\infty>0$, $m_\infty=0$, we have 

\begin{proposition} \label{Pr.3.1}
Suppose that \eqref{1.5} hold and $\int_{\Omega}\rho_0>\int_{\Omega}m_0$. Let $p_0\in(\frac 32, 3)$, $q_0\in(3,\frac{3p_0}{3-p_0})$.
There exists
$\varepsilon>0$ such that for any initial data $(\rho_0,m_0,c_0,u_0)$ fulfilling
\eqref{1.6} as well as
{\setlength\abovedisplayskip{4pt}
\setlength\belowdisplayskip{4pt}
\begin{align*}
\|\rho_0-\rho_\infty\|_{L^{p_0}(\Omega)}<\varepsilon,\quad  \|m_0\|_{L^{q_0}(\Omega)}<\varepsilon,
\quad\| c_0\|_{L^{\infty}(\Omega)}<\varepsilon, \quad\|u_0\|_{L^{3}(\Omega)}<\varepsilon,
\end{align*}}
\eqref{1.4} admits a global classical solution $(\rho,m,c,u,P)$.
In particular, 
for any $\alpha_1\in(0,\min\{\lambda_1,\rho_\infty\})$, $\alpha_2\in(0,\min\{\alpha_1,\lambda_1',1\})$, 
there exist constants $K_i$, $i=1,2,3,4$,  such that for all $t\geq 1 $
{\setlength\abovedisplayskip{4pt}
\setlength\belowdisplayskip{4pt}\begin{align}
\|m(\cdot,t)\|_{L^\infty(\Omega)}\leq & K_1e^{-\alpha_1 t},\label{3.48}
\\
\|\rho(\cdot,t)-\rho_\infty\|_{L^\infty(\Omega)}\leq & K_2e^{-\alpha_1 t},\label{3.49}
\\
\|c(\cdot,t)\|_{W^{1,\infty}(\Omega)}\leq & K_3e^{-\alpha_2t},\label{3.50}
\\
\|u(\cdot,t)\|_{L^\infty(\Omega)}\leq & K_4 e^{-\alpha_2 t}.\label{3.51}
\end{align}}
\end{proposition}
Proposition 3.1 is the consequence of the following lemmas. In the proof thereof the constants $C_i>0$, $i=1,\ldots,10$, refer to those in
Lemma 2.1--2.4, respectively.  The following verifiable observations will warrant the  choice in these lemmas. \vspace{-0.5em}
\begin{lemma}\label{Lemma3.11a}
Under the assumptions of Proposition 3.1 and $\sigma=\int_0^\infty\left(1+s^{-\frac{3}{2p_0}}\right)e^{-\alpha_1s}ds$, there exist $M_1>0,M_2>0$ and $\varepsilon\in (0,1)$ such that
\begin{align}
&C_2+2 C_2C_{10} e^{(1+C_1+C_1|\Omega|^{\frac1{p_0}-\frac1{q_0}})\sigma} \leq \displaystyle \frac {M_2}4,
\label{33.1}
\\
&C_4 C_{10} C_SM_2 (e^{(1+C_1+C_1|\Omega|^{\frac1{p_0}-\frac1{q_0}})\sigma}+\rho_\infty|\Omega|^{\frac{1}{q_0}})\leq \displaystyle \frac {M_1}8,\label{33.3}
\\
& C_6\!+\! 2C_6C_9 C_{10} \|\nabla\phi\|_{L^\infty(\Omega)}(M_1\!+\!C_1\!+\!C_1|\Omega|^{\frac1{p_0}\!-\!\frac1{q_0}}\!+4e^{(1\!+\!C_1\!+
\!C_1|\Omega|^{\frac1{p_0}\!-\!\frac1{q_0}})\sigma})
<\frac {M_3}4,\label{33.3a}
\\
& C_7\!+\! 2C_7C_9 C_{10} \|\nabla\phi\|_{L^\infty(\Omega)}|\Omega|^{\frac1{3}\!-\!\frac1{q_0}}(M_1\!+\!C_1\!+\!C_1|\Omega|^{\frac1{p_0}\!-\!\frac1{q_0}}\!+4e^{(1\!+\!C_1\!+
\!C_1|\Omega|^{\frac1{p_0}\!-\!\frac1{q_0}})\sigma})
<\frac {M_4}4,\label{33.3b}
\\
&3 C_{10} C_4C_S(M_1+C_1+C_1|\Omega|^{\frac1{p_0}-\frac1{q_0}})M_2\varepsilon\leq \displaystyle \frac {M_1}8,\label{33.4}
\\
&3 C_{10}C_4 (M_1+C_1+C_1|\Omega|^{\frac1{p_0}-\frac1{q_0}})M_3
\varepsilon
\leq \displaystyle\frac{M_1}{4},\label{33.5}
\\
&12 C_2 C_{10} M_3
\varepsilon<1,\label{33.2}
\\
&12 C_7C_9 C_{10}M_3 \varepsilon \leq 1, \label{33.5b}
\\
&12 C_6C_9 C_{10}M_4 \varepsilon \leq 1 \label{33.5a}.
\end{align}
\end{lemma}

Let
\begin{align}
\!\!T\!\triangleq\!\sup\!\left\{
\!\widetilde{T}\!\in\!(0,T_{max})\!\left|
\begin{array}{l}
\!\|(\rho\!-\!m)(\cdot,t)\!-\!e^{t\Delta}(\rho_0\!-\!m_0)\|_{L^{\theta}(\Omega)}\!\leq\! M_1\varepsilon(1\!+\!t^{-\frac{3}{2}
(\frac{1}{p_0}-\frac{1}{\theta})})e^{-\alpha_1t}\\
\hbox{for all} \,\,\theta\in[q_0,\infty], ~t\in[0,\widetilde{T});
\\
 \!\|\nabla c(\cdot,t)\|_{L^\infty(\Omega)}\leq M_2\varepsilon(1+t^{-\frac12})e^{-\alpha_1t}~\hbox{for all}~t\in[0,\widetilde{T});\\
 \|u(\cdot,t)\|_{L^{q_0}(\Omega)}\leq M_3\varepsilon(1+t^{-\frac12+\frac{3}{2q_0}}) e^{-\alpha_2 t}
~\hbox{for all}~t\in[0,\widetilde{T});\\
\|\nabla u(\cdot,t)\|_{L^{3}(\Omega)}\leq M_4\varepsilon(1+t^{-\frac12}) e^{-\alpha_2 t}
~\hbox{for all}~t\in[0,\widetilde{T}).
\end{array}
\!\!\right.
\right\}.\label{3.52}
\end{align}
Then $T>0$ is well-defined by Lemma \ref{lemma2.5} and (1.6). Now we claim that $T=T_{max}=\infty$ if $\varepsilon$ is sufficiently small.
To this end, by the contradiction argument, it only needs to verify that all of the estimates
mentioned in \eqref{3.52} also hold with even smaller coefficients on the right-side thereof, which mainly rely on $L^p-L^q$ estimates for the Neumann heat semigroup
and the fact that the classical solution on $(0,T_{max})$ can be written as 
\begin{align}
(\rho-m)(\cdot,t)=& e^{t\Delta}(\rho_0-m_0)-\!\int_0^te^{(t-s)\Delta}(\nabla\cdot(\rho \mathcal{S}(x,\rho,c)\nabla c)+ u\cdot\nabla(\rho-m))(\cdot,s)ds,\label{3.53}
\\
m(\cdot,t)=& e^{t\Delta}m_0-\int_0^te^{(t-s)\Delta}(\rho m+u\cdot\nabla m)(\cdot,s)ds,\label{3.54}
\\
c(\cdot,t)=& e^{t(\Delta-1)}c_0+\int_0^te^{(t-s)(\Delta-1)}(m-u\cdot\nabla c)(\cdot,s)ds,\label{3.55}
\\
u(\cdot,t)=& e^{-tA}u_0+\int_0^te^{-(t-s)A}\mathcal{P}((\rho+m)\nabla\phi-
 (u\cdot\nabla)u )(\cdot,s)ds \label{3.56}
\end{align}
for all $t\in(0,T_{max})$ according to the variation-of-constants formula.

Although the proof of Lemma \ref{Lemma3.2} and Lemma \ref{Lemma3.12} below is very similar to that of Lemma 3.11 and Lemma 3.12 in \cite{LPW}, respectively,
 we  give their  proofs for the convenience of the interested reader.

\begin{lemma}\label{Lemma3.2}
Under the assumptions of Proposition 3.1, for all $t\in(0,T)$ and $\theta\in[q_0,\infty]$, there exists constant $ M_5>0$ such that
\begin{align*}
\|(\rho-m)(\cdot,t)-\rho_\infty\|_{L^\theta(\Omega)}\leq M_5\varepsilon (1+t^{-\frac{3}{2}(\frac{1}{p_0}-\frac{1}{\theta})})e^{-\alpha_1t}.
\end{align*}
\end{lemma}

\proof Due to $e^{t\Delta}\rho_\infty=\rho_\infty$ and $\int_\Omega (\rho_0-m_0-\rho_\infty)=0$, the definition of $T$ and
Lemma \ref{Lemma 2.1}(i) 
show that for all $t\in(0,T)$ and $\theta\in[q_0,\infty]$,
{\setlength\abovedisplayskip{5pt}
\setlength\belowdisplayskip{4pt}
\begin{align*}
&\|(\rho-m)(\cdot,t)-\rho_\infty\|_{L^\theta(\Omega)}
\\
\leq&\|(\rho-m)(\cdot,t)-e^{t\Delta}(\rho_0-m_0)\|_{L^\theta(\Omega)}
+\|e^{t\Delta}(\rho_0-m_0-\rho_\infty)\|_{L^\theta(\Omega)}
\\
\leq& M_1\varepsilon (1+t^{-\frac{3}{2}(\frac{1}{p_0}-\frac{1}{\theta})})e^{-\alpha_1t}
+C_1(1+t^{-\frac{3}{2}(\frac{1}{p_0}-\frac{1}{\theta})})(\|\rho_0-\rho_\infty\|_{L^{p_0}(\Omega)}
+\|m_0\|_{L^{p_0}(\Omega)})e^{-\lambda_1t}
\\
\leq&M_5\varepsilon(1+t^{-\frac{3}{2}(\frac{1}{p_0}-\frac{1}{\theta})})e^{-\alpha_1t},
\end{align*}}
where $M_5=M_1+C_1+C_1|\Omega|^{\frac1{p_0}-\frac1{q_0}}$.
\begin{lemma}\label{Lemma3.12}
Under the assumptions of Proposition 3.1, for any $k>1$, 
{\setlength\abovedisplayskip{4pt}
\setlength\belowdisplayskip{4pt}
\begin{align}\label{3.57}
\|m(\cdot,t)\|_{L^k(\Omega)}\leq M_6 \|m_0\|_{L^k(\Omega)} e^{-\rho_\infty t}\quad\hbox{for all}~ t\in(0,T)
\end{align}}
with $\sigma=\int_0^\infty(1+s^{-\frac{3}{2p_0}})e^{-\alpha_1s}ds$ and $M_6=e^{M_5\sigma \varepsilon }$.
\end{lemma}\vspace{-1em}
\proof  Testing the first equation in \eqref{1.1} with $m^{k-1}$ ($k>1$) and integrating by parts, it holds that
$$
\frac{d}{dt}\int_\Omega m^k\leq-k\int_\Omega\rho m^k ~~\hbox{on}~~(0,T).
$$
In view of  $-\rho\leq|\rho-m-\rho_\infty|-m-\rho_\infty\leq-\rho_\infty+|\rho-m-\rho_\infty|$,
Lemma \ref{Lemma3.2} yields
{\setlength\abovedisplayskip{4pt}
\setlength\belowdisplayskip{4pt}\begin{align*}
\frac{d}{dt}\int_\Omega m^k&\leq-k\rho_\infty\int_\Omega m^k+k\int_\Omega m^k|\rho-m-\rho_\infty|
\\
&\leq-k\rho_\infty\int_\Omega m^k+k\|\rho-m-\rho_\infty\|_{L^\infty(\Omega)}\int_\Omega m^k
\\
&\leq-k\rho_\infty\int_\Omega m^k+kM_5\varepsilon(1+t^{-\frac{3}{2p_0}})e^{-\alpha_1t}\int_\Omega m^k
\end{align*}}
and thus
 $$
\int_\Omega m^k \leq\int_\Omega m_0^k \exp\{-k\rho_\infty t+kM_5\varepsilon \int_0^t(1+s^{-\frac{3}{2p_0}})e^{-\alpha_1s}ds\}
\leq \|m_0\|_{L^k(\Omega)} ^ke^{k(M_5\sigma\varepsilon -\rho_\infty t)},
$$
from which  \eqref{3.57} follows immediately.\vspace{-0.5em}

\begin{lemma}\label{Lemma 3.14}
Under the assumptions of Proposition 3.1, we have
$$
\|u(\cdot,t)\|_{L^{q_0}(\Omega)}\leq \frac {M_3}2 \varepsilon\left(1+t^{-\frac12+\frac{3}{2q_0}}\right) e^{-\alpha_2 t}~~\hbox{for all}~~t\in(0,T).
$$
\end{lemma}

\proof For $\alpha_2<\lambda_1'$, we fix $ \mu\in (\alpha_2, \lambda_1')$.  According to \eqref{3.56}, Lemma 2.2(ii) and Lemma 2.3,
we infer that
{\setlength\abovedisplayskip{4pt}
\setlength\belowdisplayskip{4pt}
\begin{align}
&\|u(\cdot,t)\|_{L^{q_0}(\Omega)}\nonumber
\\
\leq& C_6t^{-\frac 32\left(\frac13-\frac{1}{q_0}\right)}e^{-\mu t}\|u_0\|_{L^3(\Omega)}
+\int_0^t\|e^{-(t-s)A}\mathcal{P}((\rho+m)\nabla\phi- (u \cdot\nabla)u)(\cdot,s)\|_{L^{q_0}(\Omega)}ds\nonumber
\\
\leq& C_6t^{-\frac 32\left(\frac 13-\frac{1}{q_0}\right)}e^{-\mu t}\|u_0\|_{L^3(\Omega)}
+C_6\int_0^te^{-\mu(t-s)}\|\mathcal{P}((\rho+m-\overline{\rho+m})\nabla\phi)(\cdot,s)\|_{L^{q_0}(\Omega)}ds\nonumber\\
& +C_6\int_0^te^{-\mu(t-s)}\|\mathcal{P}((u \cdot\nabla)u)(\cdot,s)\|_{L^{q_0}(\Omega)}ds\label{3.58}
\\
\leq& C_6 t^{-\frac 12+\frac{3}{2q_0}}e^{-\mu t}\|u_0\|_{L^3(\Omega)}
+C_6C_9\|\nabla\phi\|_{L^\infty(\Omega)}\int_0^te^{-\mu(t-s)}\|(\rho+m-\overline{\rho+m})(\cdot,s)\|_{L^{q_0}(\Omega)}ds\nonumber\\
& +C_6C_9 \int_0^t(t-s)^{-\frac 12}e^{-\mu(t-s)}\|(u \cdot\nabla)u(\cdot,s)\|_{L^{\frac 1{\frac 13+\frac 1{q_0}}}(\Omega)}ds \nonumber\\
=:&C_6 t^{-\frac 12+\frac{3}{2q_0}}e^{-\mu t}\|u_0\|_{L^3(\Omega)}+J_1+J_2,\nonumber
\end{align}}
where $\mathcal{P}(\overline{\rho+m}\nabla\phi)=\overline{\rho+m} \mathcal{P}(\nabla\phi)=0$ is used.

Due to  $\alpha_1<\rho_\infty$, the application of  Lemma  \ref{Lemma3.2} and Lemma \ref{Lemma3.12} shows  that
{\setlength\abovedisplayskip{4pt}
\setlength\belowdisplayskip{4pt}
\begin{align}
J_1\leq &C_6C_9\|\nabla\phi\|_{L^\infty(\Omega)}\int_0^te^{-\mu(t-s)}\|(\rho-m-\overline{\rho-m})(\cdot,s)+2(m-\overline{m})(\cdot,s)\|_{L^{q_0}(\Omega)}ds
\nonumber\\
\leq&
C_6C_9\|\nabla\phi\|_{L^\infty(\Omega)}\int_0^te^{-\mu(t-s)}(\|(\rho-m-\rho_\infty)(\cdot,s)\|_{L^{q_0}(\Omega)}+2\|(m-\overline{m})(\cdot,s)\|_{L^{q_0}
(\Omega)})ds \nonumber\\
\leq&  C_6C_9\|\nabla\phi\|_{L^\infty(\Omega)} M_7'
\varepsilon   \int_0^te^{-\mu(t-s)} (1+s^{-\frac{3}{2}(\frac{1}{p_0}-\frac{1}{q_0})})e^{-\alpha_1s}ds
\label{1.50a}
\end{align}
with $M_7'=M_5+4e^{M_5 \sigma \varepsilon}$.

On the other hand, by the H\"{o}lder inequality and definition of $T$, we have
\begin{align}
J_2\leq &C_6C_9 \int_0^t(t-s)^{-\frac 12}e^{-\mu(t-s)}\|u(\cdot,s)\|_{L^{q_0}(\Omega)}
\|\nabla u(\cdot,s)\|_{L^{3}(\Omega)}ds
\nonumber\\
\leq&
3 C_6C_9 M_3 M_4 \varepsilon^2  \int_0^t(t-s)^{-\frac 12} e^{-\mu(t-s)}
(1+s^{-1+\frac{3}{2q_0}}) e^{-2\alpha_2 s} ds.
\label{1.50}
\end{align}
Now, plugging \eqref{1.50a}, \eqref{1.50} into \eqref{3.58} and applying Lemma 2.4, we end up with
{\setlength\abovedisplayskip{4pt}
\setlength\belowdisplayskip{4pt}
\begin{align*}
&\|u(\cdot,t)\|_{L^{q_0}(\Omega)}
\\
\leq& C_6 t^{-\frac12+\frac{3}{2q_0}}e^{-\mu t}\|u_0\|_{L^3(\Omega)}
+C_6C_9 C_{10} \|\nabla\phi\|_{L^\infty(\Omega)}M_7'\varepsilon(1+t^{\min\{0,1-\frac 32(\frac{1}{p_0}-\frac{1}{q_0})\}})e^{-\alpha_2t}\\
& +3 C_6C_9 C_{10} M_3 M_4 \varepsilon^2(1+t^{-\frac12+\frac{3}{2q_0}} )e^{-\alpha_2 t}
\\
\leq& C_6 t^{-\frac12+\frac{3}{2q_0}}e^{-\mu t}\varepsilon
+2 C_6C_9 C_{10} \|\nabla\phi\|_{L^\infty(\Omega)}M_7'\varepsilon e^{-\alpha_2t}+3 C_6C_9 C_{10} M_3 M_4 \varepsilon^2(1+t^{-\frac12+\frac{3}{2q_0}} )e^{-\alpha_2 t}
\\
\leq & \frac {M_3}2 \varepsilon(1+t^{-\frac12+\frac{3}{2q_0}}) e^{-\alpha_2 t},
\end{align*}}
where \eqref{33.3a}, \eqref{33.5a} and  the fact that $\frac{3}{2}(\frac{1}{p_0}-\frac{1}{q_0})<1$ are used.

The estimate for the gradient is also preserved.
\begin{lemma}\label{Lemma 3.14a}
Under the assumptions of Proposition 3.1, we have
$$
\|\nabla u(\cdot,t)\|_{L^{3}(\Omega)}\leq \frac {M_4}2 \varepsilon(1+t^{-\frac12})  e^{-\alpha_2 t}~~\hbox{for all}~~t\in(0,T).
$$
\end{lemma}
\proof According to \eqref{3.56}, we have
$$ \nabla u(\cdot,t)= \nabla e^{-tA}u_0+\int_0^t \nabla e^{-(t-s)A}(\mathcal{P}((\rho+m)\nabla\phi)-
\mathcal{P}( (u\cdot\nabla)u ))(\cdot,s)ds. $$
Applying Lemma 2.2(iii), Lemma 2.3 and the H\"{o}lder inequality, we arrive at
{\setlength\abovedisplayskip{4pt}
\setlength\belowdisplayskip{4pt}
\begin{align}
&\| \nabla u(\cdot,t)\|_{L^{3}(\Omega)}\nonumber
\\
\leq& C_7 t^{-\frac 12}e^{-\mu t}\|u_0\|_{L^3(\Omega)}
+\int_0^t\| \nabla e^{-(t-s)A}\mathcal{P}((\rho+m)\nabla\phi- (u \cdot\nabla)u)(\cdot,s)\|_{L^3(\Omega)}ds\nonumber
\\
\leq& C_7 t^{-\frac 12}e^{-\mu t}\varepsilon
+C_7\int_0^t (t-s)^{-\frac{1}{2}}e^{-\mu(t-s)}\|\mathcal{P}((\rho+m-\overline{\rho+m})\nabla\phi)(\cdot,s)\|_{L^{3}(\Omega)}ds\nonumber\\
& +C_7\int_0^t  (t-s)^{-\frac 12-\frac{3}{2q_0}}  e^{-\mu(t-s)}\|\mathcal{P}((u \cdot\nabla)u)(\cdot,s)\|_
{L^{\frac {3q_0}{3+q_0}}(\Omega)
}
ds\label{3.58d}
\\
\leq&  C_7 t^{-\frac 12}e^{-\mu t}\varepsilon
+C_7C_9\|\nabla\phi\|_{L^\infty(\Omega)} |\Omega|^{\frac1{3}\!-\!\frac1{q_0}}\int_0^t  (t-s)^{-\frac{1}{2}} e^{-\mu(t-s)}\|(\rho+m-\overline{\rho+m})(\cdot,s)\|_{L^{q_0}(\Omega)}ds\nonumber\\
& +C_7C_9\int_0^t  (t-s)^{-\frac 12-\frac{3}{2q_0}}  e^{-\mu(t-s)} \|\nabla u(\cdot,s)\|_{L^{3}(\Omega)}\|u(\cdot,s)\|_{L^{q_0}(\Omega)}ds
 \nonumber\\
=:&C_7 t^{-\frac 12}e^{-\mu t}\varepsilon+\jmath_1+\jmath_2,\nonumber
\end{align}
where $\mathcal{P}(\overline{\rho+m}\nabla\phi)=\overline{\rho+m} \mathcal{P}(\nabla\phi)=0$ is used.

Due to  $\alpha_1<\rho_\infty$, the application of  Lemma  \ref{Lemma3.2} and Lemma \ref{Lemma3.12} shows  that
{\setlength\abovedisplayskip{4pt}
\setlength\belowdisplayskip{4pt}
\begin{align}
\jmath_1\leq &C_7C_9\|\nabla\phi\|_{L^\infty(\Omega)} |\Omega|^{\frac1{3}\!-\!\frac1{q_0}}\int_0^t (t-s)^{-\frac{1}{2}}e^{-\mu(t-s)}\|(\rho-m-\overline{\rho-m})(\cdot,s)+2(m-\overline{m})(\cdot,s)\|_{L^{q_0}(\Omega)}ds
\nonumber\\
\leq&
C_7C_9\|\nabla\phi\|_{L^\infty(\Omega)}  |\Omega|^{\frac1{3}\!-\!\frac1{q_0}} \int_0^t(t-s)^{-\frac{1}{2}}e^{-\mu(t-s)}(\|(\rho-m-\rho_\infty)(\cdot,s)\|_{L^{q_0}(\Omega)}+2\|(m-\overline{m})(\cdot,s)\|_{L^{q_0}
(\Omega)})ds \nonumber\\
\leq&  C_7C_9\|\nabla\phi\|_{L^\infty(\Omega)} |\Omega|^{\frac1{3}\!-\!\frac1{q_0}}M_7'
\varepsilon   \int_0^te^{-\mu(t-s)} (1+s^{-\frac{3}{2}(\frac{1}{p_0}-\frac{1}{q_0})})(t-s)^{-\frac{1}{2}}e^{-\alpha_1s}ds.
\label{1.50e}
\end{align}
On the other hand, from the H\"{o}lder inequality and definition of $T$, it follows that
\begin{align}
\jmath_2
\leq
3 C_7C_9 M_3 M_4 \varepsilon^2  \int_0^t(t-s)^{-\frac 12-\frac{3}{2q_0}}   e^{-\mu(t-s)}
(1+s^{-1+\frac{3}{2q_0}}) e^{-2\alpha_2 s} ds.
\label{1.50d}
\end{align}
Therefore, inserting \eqref{1.50d}, \eqref{1.50e} into \eqref{3.58d} and applying Lemma 2.4, we get
{\setlength\abovedisplayskip{4pt}
\setlength\belowdisplayskip{4pt}
\begin{align*}
&\|\nabla u(\cdot,t)\|_{L^{q_0}(\Omega)}
\\
\leq& C_7 t^{-\frac12}e^{-\mu t}\varepsilon
+C_7C_9 C_{10} \|\nabla\phi\|_{L^\infty(\Omega)}|\Omega|^{\frac1{3}\!-\!\frac1{q_0}}M_7'\varepsilon(1+t^{\min\{0,
\frac12-\frac{3}{2}(\frac{1}{p_0}-\frac{1}{q_0})\}})e^{-\alpha_2t}\\
& +3 C_7C_9 C_{10} M_3 M_4 \varepsilon^2(1+t^{-\frac12} )e^{-\alpha_2 t}
\\
\leq& C_7 t^{-\frac12}e^{-\mu t}\varepsilon
\!+\!2 C_7C_9 C_{10} \|\nabla\phi\|_{L^\infty(\Omega)}|\Omega|^{\frac1{3}\!-\!\frac1{q_0}}M_7' \varepsilon e^{-\alpha_2t}\!+3 C_7C_9 C_{10} M_3 M_4 \varepsilon^2(1+t^{-\frac12} )e^{-\alpha_2 t}
\\
\leq & \frac {M_4}2 \varepsilon(1+t^{-\frac12}) e^{-\alpha_2 t},
\end{align*}}
where \eqref{33.3b}, \eqref{33.5b} and  the fact that $q_0\in(3,\frac{3p_0}{3-p_0}), p_0\in (\frac{3}{2},3)$ are used.

\begin{lemma}\label{lemma3.15}
Under the assumptions of Proposition 3.1, we have
\begin{align*}
\|\nabla c(\cdot,t)\|_{L^{\infty}(\Omega)}\leq \frac{M_2}{2}\varepsilon (1+t^{-\frac12}) e^{-\alpha_1 t}~~\hbox{for all}~~t\in(0,T).
\end{align*}

\end{lemma}

\proof By \eqref{3.55} and Lemma 2.1(ii), we have
\begin{align}
\|\nabla c(\cdot,t)\|_{L^\infty(\Omega)}
&\leq\|e^{t(\Delta-1)}\nabla c_0\|_{L^\infty(\Omega)}+\int_0^t\|\nabla
e^{(t-s)(\Delta-1)}(m-u\cdot\nabla c)(\cdot,s)\|_{L^\infty(\Omega)}ds\nonumber
\\
&\leq C_2(1+t^{-\frac12})e^{-(\lambda_1+1)t}\|c_0\|_{L^\infty(\Omega)}+\int_0^t\|\nabla
e^{(t-s)(\Delta-1)}m(\cdot,s)\|_{L^\infty(\Omega)}ds\nonumber
\\
&\quad+\int_0^t\|\nabla
e^{(t-s)(\Delta-1)}u\cdot\nabla c(\cdot,s)\|_{L^\infty(\Omega)}ds.\label{3.60}
\end{align}
Now we estimate the last two integrals on the right-side of the above inequality. From Lemma 2.1(ii), Lemma 2.4,  Lemma \ref{Lemma3.12} with $k=q_0$  and the fact that $q_0>3$, it follows that{\setlength\abovedisplayskip{4pt}
\setlength\belowdisplayskip{4pt}
\begin{align}\label{3.61}
\int_0^t\|\nabla
e^{(t-s)(\Delta-1)}m\|_{L^\infty(\Omega)}ds
\leq&C_2\int_0^t(1+(t-s)^{-\frac12-\frac{3}{2q_0}})e^{-(\lambda_1+1)(t-s)}\|m\|_{L^{q_0}(\Omega)}ds
\\
\leq&C_2 M_6\varepsilon
\int_0^t(1+(t-s)^{-\frac12-\frac{3}{2q_0}})e^{-(\lambda_1+1)(t-s)}
e^{-\rho_\infty s}ds
\nonumber
\\
\leq& C_2C_{10} M_6
(1+t^{\min\{0,\frac12-\frac{3}{2q_0}\}})\varepsilon e^{-\alpha_1t}\nonumber
\\
\leq& 2 C_2C_{10} M_6 \varepsilon e^{-\alpha_1 t}.\nonumber
\end{align}}
On the other hand, by Lemma 2.1(ii), Lemma 2.4, Lemma \ref{Lemma 3.14} and the definition of $T$, we obtain{\setlength\abovedisplayskip{4pt}
\setlength\belowdisplayskip{4pt}
\begin{align}
&\int_0^t\|\nabla
e^{(t-s)(\Delta-1)}u\cdot\nabla c\|_{L^\infty(\Omega)}ds\nonumber
\\
\leq&C_2\int_0^t(1+(t-s)^{-\frac12-\frac{3}{2q_0}})e^{-(\lambda_1+1)(t-s)}\|u\cdot\nabla c\|_{L^{q_0}(\Omega)}ds\label{3.62}
\\
\leq&C_2\int_0^t(1+(t-s)^{-\frac12-\frac{3}{2q_0}})e^{-(\lambda_1+1)(t-s)}\|u\|_{L^{q_0}(\Omega)}\|\nabla c\|_{L^\infty(\Omega)}ds\nonumber
\\
\leq& C_2 M_3M_2 \varepsilon^2 \int_0^t(1+(t-s)^{-\frac12-\frac{3}{2q_0}})e^{-(\lambda_1+1)(t-s)}(1+s^{-\frac12+\frac{3}{2q_0}}) (1+s^{-\frac12}) e^{-(\alpha_1+\alpha_2) s}ds \nonumber
\\
\leq&3C_2M_3M_2\varepsilon^{2}\int_0^te^{-(\lambda_1+1)(t-s)}
 e^{-(\alpha_1+\alpha_2) s}(1+(t-s)^{-\frac12-\frac{3}{2q_0}})(1+s^{-1+\frac{3}{2q_0}})ds\nonumber
\\
\leq&3C_2 C_{10}M_2  M_3\varepsilon^2(1+t^{-\frac12})e^{-\alpha_1 t}.\nonumber
\end{align}}
From \eqref{3.60}--\eqref{3.62}, it follows that
\begin{align*}
\|\nabla c\|_{L^\infty(\Omega)}
&\leq (C_2+2 C_2C_{10} M_6 +3C_2 C_{10}M_2  M_3\varepsilon)
(1+t^{-\frac12})\varepsilon e^{-\alpha_1t}
\\
&\leq \frac{M_2}{2}(1+t^{-\frac12})\varepsilon e^{-\alpha_1 t},
\end{align*}
due to the choice of $M_2, M_3$ and $\varepsilon$ in \eqref{33.1} and \eqref{33.2}, and thereby completes the proof.

\begin{lemma}\label{Lemma3.16}
Under the assumptions of Proposition 3.1, for all $\theta\in[q_0,\infty]$ and $t\in(0,T)$,
\begin{align*}
\|(\rho-m)(\cdot,t)-e^{t\Delta}(\rho_0-m_0)\|_{L^\theta(\Omega)}\leq \frac{M_1}{2}\varepsilon (1+t^{-\frac{3}{2}(\frac{1}{p_0}-\frac{1}{\theta})}) e^{-\alpha_1 t}.
\end{align*}
\end{lemma}
\proof According to \eqref{3.53}, Lemma 2.1(iv), we have
{\setlength\abovedisplayskip{4pt}
\setlength\belowdisplayskip{4pt}\begin{align*}
&\|(\rho-m)(\cdot,t)-e^{t\Delta}(\rho_0-m_0)\|_{L^\theta(\Omega)}
\\
\leq&\int_0^t\|e^{(t-s)\Delta}(\nabla\cdot(\rho \mathcal{S}(x,\rho,c)\nabla c)+u\cdot\nabla(\rho-m))(\cdot,s)\|_{L^\theta(\Omega)}ds
\\
\leq&\int_0^t\|e^{(t-s)\Delta}\nabla\cdot(\rho\mathcal{S}(x,\rho,c)\nabla c)(\cdot,s)\|_{L^\theta(\Omega)}ds
+\int_0^t\|e^{(t-s)\Delta}\nabla\cdot((\rho-m-\rho_\infty)u)(\cdot,s)\|_{L^\theta(\Omega)}ds
\\
\leq&C_4C_S\int_0^t(1+(t-s)^{-\frac12-\frac{3}{2}(\frac{1}{q_0}-\frac{1}{\theta})})
e^{-\lambda_1(t-s)}\|\rho(\cdot,s)\|_{L^{q_0}(\Omega)}\|\nabla c(\cdot,s)\|_{L^\infty(\Omega)}ds
\\
&+C_4\int_0^t(1+(t-s)^{-\frac12-\frac{3}{2}(\frac{1}{q_0}-\frac{1}{\theta})})
e^{-\lambda_1(t-s)}\|u(\rho-m-\rho_\infty)(\cdot,s)\|_{L^{q_0}(\Omega)}ds
\\
=:&I_1+I_2.
\end{align*}}
Now we need to estimate $I_1$ and $I_2$. Firstly, from Lemma \ref{Lemma3.2} and Lemma \ref{Lemma3.12}, we obtain
\begin{align}
\|\rho(\cdot,s)\|_{L^{q_0}(\Omega)}&\leq \|(\rho-m-\rho_\infty)(\cdot,s)\|_{L^{q_0}(\Omega)}+\|m(\cdot,s)\|_{L^{q_0}(\Omega)}+\|\rho_\infty\|_{L^{q_0}(\Omega)}\label{3.63}
\\
&\leq M_5\varepsilon
(1+s^{-\frac{3}{2}\left(\frac{1}{p_0}-\frac{1}{q_0}\right)})e^{-\alpha_1s}+M_8 \nonumber
\end{align}
with $M_8=e^{(1+C_1+C_1|\Omega|^{\frac1{p_0}-\frac1{q_0}})\sigma}+\rho_\infty|\Omega|^{\frac{1}{q_0}}$,
which along with Lemma \ref{lemma3.15} and Lemma 2.1 implies that\vspace{-0.5em}
{\setlength\abovedisplayskip{4pt}
\setlength\belowdisplayskip{4pt}\begin{align}
I_1&\leq C_4C_SM_8\int_0^t(1+(t-s)^{-\frac12-\frac{3}{2}(\frac{1}{q_0}-\frac{1}{\theta})})
e^{-\lambda_1(t-s)}\|\nabla c\|_{L^\infty(\Omega)}ds
\label{1.52}\\
 &+ C_4C_SM_5\varepsilon
\int_0^t(1+(t-s)^{-\frac12-\frac{3}{2}(\frac{1}{q_0}-\frac{1}{\theta})})(1+s^{-\frac{3}{2}(\frac{1}{p_0}-\frac{1}{q_0})})
e^{-\alpha_1s}
e^{-\lambda_1(t-s)}\|\nabla c\|_{L^\infty(\Omega)}ds
 \nonumber\\
&\leq C_4C_SM_8M_2\varepsilon \int_0^t(1+(t-s)^{-\frac12-\frac{3}{2}(\frac{1}{q_0}-\frac{1}{\theta})})
e^{-\lambda_1(t-s)}(1+s^{-\frac12}) e^{-\alpha_1 s}ds \nonumber\\
& +3 C_4C_S M_5 M_2\varepsilon^2
\int_0^t(1+(t-s)^{-\frac12-\frac{3}{2}(\frac{1}{q_0}-\frac{1}{\theta})})(1+s^{-\frac{1}{2}-\frac{3}{2}(\frac{1}{p_0}-\frac{1}{q_0})})
e^{-2\alpha_1s}
e^{-\lambda_1(t-s)}ds\nonumber
\\
&\leq  C_{10}C_4C_S (M_8M_2+3 M_5M_2\varepsilon)(1+t^{-\frac{3}{2}(\frac{1}{p_0}-\frac{1}{\theta})})\varepsilon e^{-\alpha_1 t}\nonumber
\\
&\leq \frac{M_1}{4}(1+t^{-\frac{3}{2}(\frac{1}{p_0}-\frac{1}{\theta})})\varepsilon e^{-\alpha_1 t},\nonumber
\end{align}}
where we have used \eqref{33.3} and \eqref{33.4} and  $\frac{1}{p_0}-\frac{1}{q_0}<\frac{1}{3}$.

On the other hand,  from  Lemma \ref{Lemma3.2} and Lemma \ref{Lemma 3.14}, it follows that
{\setlength\abovedisplayskip{4pt}
\setlength\belowdisplayskip{4pt}\begin{align}
I_2&=C_4\int_0^t(1+(t-s)^{-\frac12-\frac{3}{2}(\frac{1}{q_0}-\frac{1}{\theta})})
e^{-\alpha_1(t-s)}\|\rho-m-\rho_\infty\|_{L^\infty(\Omega)}\|u\|_{L^{q_0}(\Omega)}ds\nonumber
\\
&\leq 3C_4M_3M_5\varepsilon^{2}\int_0^t(1+(t-s)^{-\frac12-\frac{3}{2}(\frac{1}{q_0}-\frac{1}{\theta})})
e^{-\alpha_1(t-s)}(1+s^{-\frac12+\frac{3 }{ 2q_0}-\frac{3}{2p_0}})e^{-(\alpha_1+\alpha_2)s}
ds\nonumber
\\
&\leq 3C_4C_{10}M_3M_5 \varepsilon^{2}(1+t^{\min\{0,\frac{3}{2}
(\frac{1}{\theta}-\frac{1}{p_0})\}})e^{-\alpha_1 t }\nonumber
\\
&\leq \frac{M_1}{4}\varepsilon (1+t^{-\frac{3}{2}(\frac{1}{p_0}-\frac{1}{\theta})})e^{-\alpha_1 t},\label{1.53}
\end{align}}
where we have used \eqref{33.5} and  $\frac{1}{p_0}-\frac{1}{q_0}<\frac{1}{3}$. Hence combining the above inequalities leads to our conclusion immediately.

Now we are ready to complete the proof of  Theorem 1.1 in the case $\mathcal{S}=0$ on $\partial\Omega$.

{\bf Proof of Proposition 3.1.}~
First from Lemma  \ref{Lemma 3.14}--\ref{Lemma3.16} and Definition \eqref{3.52}, it follows that $T=T_{max}$.
It remains to show  that $T_{max}=\infty$ and convergence result asserted in Proposition 3.1.
Supposed that   $T_{max}<\infty$,
we only need  to show that for all  $t \leq T_{max}$,
{\setlength\abovedisplayskip{4pt}
\setlength\belowdisplayskip{4pt}
  $$ \|\rho(\cdot,t)\|_{L^\infty(\Omega)}+\|m(\cdot,t)\|_{L^\infty(\Omega)}
+\|c(\cdot,t)\|_{W^{1,\infty}(\Omega)}+
\|A^{\beta}u(\cdot,t)\|_{L^2(\Omega)}< \infty$$
} with $\beta\in(\frac{3}{4},1)$ according to the extensibility criterion in Lemma \ref{lemma2.5}.

Let $t_0:=\min\{1,\frac{T_{max}}3\}$. Then from Lemma \ref{Lemma3.12}, there  exists  $K_1>0$ such that for $t\in(t_0,T_{max})$,
{\setlength\abovedisplayskip{4pt}
\setlength\belowdisplayskip{4pt}\begin{align}\label{3.66a}
\|m(\cdot,t)\|_{L^\infty(\Omega)}\leq K_1e^{-\rho_\infty t}.
\end{align}}
Moreover, from Lemma \ref{Lemma3.2} and the fact that
{\setlength\abovedisplayskip{4pt}
\setlength\belowdisplayskip{4pt}
$$\|\rho(\cdot,t)-\rho_\infty\|_{L^\infty(\Omega)}\leq\|(\rho-m)(\cdot,t)-\rho_\infty\|_{L^\infty(\Omega)}+\|m(\cdot,t)\|_{L^\infty(\Omega)},
$$}
it follows that for  all $t\in(t_0,T_{max})$ and  some constant $K_2>0$,
{\setlength\abovedisplayskip{4pt}
\setlength\belowdisplayskip{4pt}\begin{align}\label{3.67a}
\|\rho(\cdot,t)-\rho_\infty\|_{L^\infty(\Omega)}\leq K_2e^{-\alpha_1 t}.
\end{align}}
Furthermore, Lemma \ref{lemma3.15} implies that there exists $K_3'>0$ such that
{\setlength\abovedisplayskip{4pt}
\setlength\belowdisplayskip{4pt}\begin{align}\label{3.66b}
\|\nabla c(\cdot,t)\|_{L^{\infty}(\Omega)}\leq K_3'e^{-\alpha_1 t}\quad \hbox{for all}\,\,t\in(t_0,T_{max})
\end{align}}
Hence  it only remains to show that
 $$
 \|c(\cdot,t)\|_{L^\infty(\Omega)}+
\|A^{\beta}u(\cdot,t)\|_{L^2(\Omega)}\leq C~~\hbox{for all}~t\in(t_0,T_{max}).
$$
for some constant $C>0$.
In fact,  we will show that
{\setlength\abovedisplayskip{4pt}
\setlength\belowdisplayskip{4pt}\begin{align}\label{3.66}
\|A^\beta u(\cdot,t)\|_{L^2(\Omega)}\leq C e^{-\alpha_2 t}
\end{align}}
for $t_0<t<T_{max}$ with some constant $C>0$.

By \eqref{3.56}, we have
{\setlength\abovedisplayskip{4pt}
\setlength\belowdisplayskip{4pt}\begin{align}
\|A^\beta u(\cdot,t)\|_{L^2(\Omega)}\leq &\|A^\beta e^{-tA} u_0\|_{L^2(\Omega)}
+\int_{0}^t\|A^\beta e^{-(t-s)A}\mathcal{P}((\rho+m-\rho_\infty)\nabla\phi)(\cdot,s)\|_{L^2(\Omega)}ds\nonumber\\
&+\int_{0}^t\|A^\beta e^{-(t-s)A}\mathcal{P}((u\cdot\nabla)u )(\cdot,s)\|_{L^2(\Omega)}ds\label{3.66a}.
\end{align}}
  According to Lemma 2.2,
$$
\|A^\beta e^{-tA} u_0\|_{L^2(\Omega)}\leq  C_5 e^{-\mu t}\|A^\beta u_0\|_{L^2(\Omega)}~~\hbox{for all}~~t\in(0,T_{max}).
$$
From Lemma 2.2, 2.3,  
 \ref{Lemma3.2} and the H\"{o}lder inequality, it follows that there exists $l_1>0$ such that 
\begin{align*}
&\int_0^t\|A^\beta e^{-(t-s)A}\mathcal{P}((\rho+m-\rho_\infty)\nabla\phi)(\cdot,s)\|_{L^2(\Omega)}ds
\\
\leq&C_5 C_9\|\nabla\phi\|_{L^\infty(\Omega)}|\Omega|^{\frac{q_0-2}{2q_0}}
\int_0^t (\|(\rho-m-\rho_\infty)(\cdot,s)\|_{L^{q_0}(\Omega)}+2\|m(\cdot,s)\|_{L^{q_0}(\Omega)})(t-s)^{-\beta}e^{-\mu(t-s)}ds
\\
\leq& C_5 C_9\|\nabla\phi\|_{L^\infty(\Omega)}|\Omega|^{\frac{q_0-2}{2q_0}}l_1
\int_0^t e^{-\mu(t-s)}(t-s)^{-\beta} (1+s^{-\frac{3}{2}(\frac{1}{p_0}-\frac{1}{q_0})})e^{-\alpha_1s}
ds\\
\leq& C_5 C_9 C_{10}\|\nabla\phi\|_{L^\infty(\Omega)}|\Omega|^{\frac{q_0-2}{2q_0}} l_1
e^{-\alpha_2 t}(1+t^{\min\{0,1-\beta-\frac{3}{2}(\frac{1}{p_0}-\frac{1}{q_0})\}}).
\end{align*}
On the other hand, let $M(t):=e^{-\alpha_2 t}\|A^\beta u(\cdot,t)\|_{L^2(\Omega)} $ for $ 0<t<T_{max}$.
By  Lemma 2.2(iv) and the Gagliardo--Nirenberg type inequality, one can see that
 \begin{align*}
 \|(u\cdot\nabla)u(\cdot,s)\|_{L^{2}(\Omega)}\leq
 & |\Omega|^{\frac 16}\|u(\cdot,s)\|_{L^{\infty}(\Omega)}\|\nabla u(\cdot,s)\|_{L^{3}(\Omega)}\\
 \leq  &l_2 \|A^\beta u(\cdot,s)\|_{L^2(\Omega)}^{\vartheta}\|u(\cdot,s)\|^{1-\vartheta}_{L^{q_0}(\Omega)} \|\nabla u(\cdot,s)\|_{L^{3}(\Omega)}
\end{align*}
for some $l_2>0$ with $\vartheta=\frac 1 {q_0}/(\frac 1 {q_0}-\frac 12+\frac {2\beta} {3})$, and thereby the application of Lemma
   2.2, \ref{Lemma2.3}, \ref{Lemma 3.14} and \ref{Lemma 3.14a} gives
 \begin{align*}
&\int_0^t\|A^\beta e^{-(t-s)A}\mathcal{P}((u\cdot\nabla)u)(\cdot,s)\|_{L^2(\Omega)}ds
\\
\leq&C_5 C_9 l_2
\int_0^t\|A^\beta u(\cdot,s)\|_{L^2(\Omega)}^{\vartheta}\|u(\cdot,s)\|^{1-\vartheta}_{L^{q_0}(\Omega)} \|\nabla u(\cdot,s)\|_{L^{3}(\Omega)}
\\
\leq& l_3 (\max_{0\leq s<T_{max}}M(s))^\vartheta
\int_0^t e^{-\mu(t-s)}(t-s)^{-\beta} (1+s^{-\frac 12+(-\frac 12+\frac 3{2q_0})(1-\vartheta)})e^{-2\alpha_2s}
ds\\
\leq&  C_{10}  l_3 (\max_{0\leq s<T_{max}}M(s))^\vartheta
(1+t^{\min\{0,\frac 12-\beta+(\frac 3{2q_0}-\frac 12)(1-\vartheta)\}})e^{-\alpha_2 t}
\end{align*}
for some $l_3>0$.

Hence inserting the above inequalities into \eqref{3.66a}, we arrive at
\begin{align*}
M(t)\leq& C_5 \|A^\beta u_0\|_{L^2(\Omega)}+
C_5 C_9 C_{10}\|\nabla\phi\|_{L^\infty(\Omega)}|\Omega|^{\frac{q_0-2}{2q_0}} l_1
(1+t^{\min\{0,1-\beta-\frac{3}{2}(\frac{1}{p_0}-\frac{1}{q_0})\}})\\
& + C_{10}  l_3 (\max_{0\leq s<T_{max}}M(s))^\vartheta
(1+t^{\min\{0,\frac 12-\beta+(\frac 3{2q_0}-\frac 12)(1-\vartheta)\}}),
\end{align*}
which implies that for some $l_4>0$ depending on $t_0$,  we have
\begin{align*}
\max_{t_0\leq t<T_{max}}M(t)\leq l_4+l_4 (\max_{0\leq t<T_{max}}M(t))^\vartheta.
\end{align*}
On the other hand, from  Lemma \ref{lemma2.5},
$
\displaystyle\max_{0\leq t\leq  t_0}M(t)\leq  l_5.
$
Therefore, we get
\begin{align*}
\max_{0\leq t<T_{max}}M(t)\leq l_4+l_5+l_4 (\max_{0\leq t<T_{max}}M(t))^\vartheta.
\end{align*}
Due to $\vartheta<1$, we infer that $M(t)\leq l_6$ for all  $t\in (0,T_{max})$ for some $l_6>0$ independent of $T_{max}$
 hence arrive at
\eqref{3.66}.

 Furthermore, due to $D(A^\beta)\hookrightarrow L^\infty(\Omega)$ with $\beta\in(\frac{3}{4},1)$ and Lemma 3.4, we get
{\setlength\abovedisplayskip{4pt}
\setlength\belowdisplayskip{4pt}
\begin{align}\label{3.67}
\|u(\cdot,t)\|_{L^\infty(\Omega)}\leq K_4  e^{-\alpha_2 t} \quad\hbox{for some}\,K_4>0 \,\hbox{and }\,  t\in(0,T_{max}).
\end{align}}
Now we turn to show that
 there exists $K_3''>0$ such that
 {\setlength\abovedisplayskip{4pt}
\setlength\belowdisplayskip{4pt}
\begin{align}\label{3.71}
\|c(\cdot,t)\|_{L^{\infty}(\Omega)}\leq K_3'' e^{-\alpha_2t}\quad\hbox{for all}\,\, t\in(0,T_{max}).
\end{align}}
From \eqref{3.55}, it follows that
\begin{align}
\|c\|_{L^\infty(\Omega)}
&\leq\|e^{t(\Delta-1)}c_0\|_{L^\infty(\Omega)}+\int_0^t\|
e^{(t-s)(\Delta-1)}(m-u\cdot\nabla c)\|_{L^\infty(\Omega)}ds\nonumber
\\
&\leq  e^{-t}\|c_0\|_{L^\infty(\Omega)}+\int_0^t\|
e^{(t-s)(\Delta-1)}m(\cdot,s)\|_{L^\infty(\Omega)}ds\label{3.72}
\\
&\quad+\int_0^t\|
e^{(t-s)(\Delta-1)}u\cdot\nabla c(\cdot,s)\|_{L^\infty(\Omega)}ds.\nonumber
\end{align}
 An application of \eqref{3.57} with $k=\infty$ yields
{\setlength\abovedisplayskip{4pt}
\setlength\belowdisplayskip{4pt}\begin{align}
\int_0^t\|e^{(t-s)(\Delta-1)}m(\cdot,s)\|_{L^\infty(\Omega)}ds
&\leq \int_0^t e^{-(t-s)}\|m(\cdot,s)\|_{L^{\infty}(\Omega)}ds\label{3.73}
\\
&\leq \|m_0\|_{L^\infty(\Omega)}M_6\int_0^t e^{-(t-s)} e^{-\rho_\infty s}
ds\nonumber
\\
& \leq \|m_0\|_{L^\infty(\Omega)}M_6C_{10} e^{-\alpha_2t}.\nonumber
\end{align}}
On the other hand,
 from  \eqref{3.67} and \eqref{3.66b}, we can see  that
 {\setlength\abovedisplayskip{4pt}
\setlength\belowdisplayskip{4pt}
\begin{align}\label{3.74}
\int_0^t\|
e^{(t-s)(\Delta-1)} u\cdot\nabla c\|_{L^\infty(\Omega)}ds
&\leq \int_0^t e^{-(t-s)}\|u\|_{L^{\infty}(\Omega)}\|\nabla c\|_{L^\infty(\Omega)}ds
\\
&\leq K_3' K_4  \int_0^t e^{-(\alpha_1+\alpha_2) s} e^{-(t-s)}ds
\nonumber
\\ &\leq K_3' K_4 C_{10} e^{-\alpha_2 t}.
 \nonumber
\end{align}}
Hence, inserting  \eqref{3.73}, \eqref{3.74} into  \eqref{3.72}, we arrive at the conclusion \eqref{3.71}. Therefore
we have $T_{max}=\infty$, and the decay estimates in \eqref{3.48}--\eqref{3.51} follow from (3.32)--(3.35) and (3.38), respectively. 


As for the case $\int_{\Omega}\rho_0<\int_{\Omega}m_0$, i.e.,
$m_\infty>0$, $\rho_\infty=0$, we also have
\begin{proposition}
Assume that \eqref{1.5} and $\int_{\Omega}\rho_0<\int_{\Omega}m_0$ hold, and let  $p_0\in(2,3)$, $q_0\in(3,\frac{3p_0}{2(3-p_0)})$. 
Then there exists
$\varepsilon>0$ such that for any initial data $(\rho_0,m_0,c_0,u_0)$ fulfilling
\eqref{1.6} as well as
{\setlength\abovedisplayskip{4pt}
\setlength\belowdisplayskip{4pt}
$$\|\rho_0\|_{L^{p_0}(\Omega)}\leq\varepsilon,\quad  \|m_0-m_\infty\|_{L^{q_0}(\Omega)}\leq\varepsilon,
\quad\|\nabla c_0\|_{L^{3}(\Omega)}\leq\varepsilon, \quad\|u_0\|_{L^{3}(\Omega)}\leq\varepsilon,
$$}
\eqref{1.4} admits a global classical solution $(\rho,m,c,u,P)$.
Furthermore,  for any $\alpha_1\!\in\!(0,\min\{\lambda_1,m_\infty,1\})$, $\alpha_2\!\in\!(0,\min\{\alpha_1,\lambda_1'\})$, there exist constants $K_i>0$, $i=1,2,3,4$, such that
\begin{align}
&\|m(\cdot,t)-m_\infty\|_{L^\infty(\Omega)}\leq K_1e^{-\alpha_1 t},\label{3.75}
\\
&\|\rho(\cdot,t)\|_{L^\infty(\Omega)}\leq K_2e^{-\alpha_1 t},\label{3.76}
\\
&\|c(\cdot,t)-m_\infty\|_{W^{1,\infty}(\Omega)}\leq K_3e^{-\alpha_2t},\label{3.77}
\\
&\|u(\cdot,t)\|_{L^\infty(\Omega)}\leq K_4 e^{-\alpha_2 t}.\label{3.78}
\end{align}
\end{proposition}

The  basic strategy in the proof of Proposition 3.2  parallels that  in the proof of Proposition 3.1 to a certain extent.  However, due to differences in the properties of $\rho$ and $m$,
there are significant differences in the details of their proofs. Thus for the convenience of the reader, we will sketch the  proof of Proposition 3.2.

The following elementary observations can be also verified easily:\vspace{-1em}
\begin{lemma}\label{Lemma3.17a}
Under the assumptions of Proposition 3.2, 
 it is possible to choose $M_1>0,M_2>0$ and $\varepsilon>0$ such that
\begin{align}
&C_3+
C_2 C_{10}(1+ C_1+C_1|\Omega|^{\frac1{p_0}-\frac1{q_0}}+M_1)\leq \displaystyle \frac {M_2}4,\label{34.1}
\\
&C_6 +2C_6C_9 C_{10} (M_1\!+2+2C_1\!+2 C_1|\Omega|^{\frac1{p_0}\!-\!\frac1{q_0}}\!)
\|\nabla\phi\|_{L^\infty(\Omega)}< \displaystyle\frac {M_3}4 \label{34.1b}
\\
&C_7 +2C_7C_9 C_{10} (M_1\!+2+2C_1\!+ 2C_1|\Omega|^{\frac1{p_0}\!-\!\frac1{q_0}}\!)
\|\nabla\phi\|_{L^\infty(\Omega)}|\Omega|^{\frac1{3}\!-\!\frac1{q_0}}< \displaystyle\frac {M_4}4 \label{34.1c}
\\
&12C_2C_{10}M_3\varepsilon\leq 1, \label{34.2}
\\
&2C_1+(\min\{1,|\Omega|\})^{-\frac 1{p_0}}\leq \displaystyle\frac {M_1}8,\quad
12C_6C_9C_{10}M_4 \varepsilon<1,\label{34.3}
\\
& 24C_4C_SC_{10}M_2\varepsilon<1,\label{34.3d}\\
& 12C_7C_9C_{10}M_3 \varepsilon<1,\label{34.3b}
\\
&12C_4C_{10} C_S M_1M_2 \varepsilon<1,\label{34.4}\\
&24C_1C_{10}(1+ C_1+C_1|\Omega|^{\frac1{p_0}-\frac1{q_0}}+M_1)\varepsilon<1,\label{34.5a}
\\
&18C_4C_{10} M_3 \varepsilon< 1.\label{34.6a}
\\
&12 C_{10} C_4M_3(1\!+\! C_1\!+\!C_1|\Omega|^{\frac1{p_0}\!-\!\frac1{q_0}})
\varepsilon \!< 1.\label{34.7}
\end{align}
\end{lemma}

Define{\setlength\abovedisplayskip{4pt}
\begin{align}\label{3.79}
T\!:=\!\sup\!\left\{\!\widetilde{T}\!\in\!(0,T_{max})\!\left|
\begin{array}{ll}
\!\|(m\!-\!\rho)(\cdot,t)\!-\!e^{t\Delta}(m_0\!-\!\rho_0)\|_{L^{\theta}(\Omega)}\!\leq\! \varepsilon(1+t^{-\frac{ 3}{2}(\frac{1}{p_0}-\frac{1}{\theta})})e^{-\alpha_1t};
\\
\!\|\rho(\cdot,t)\|_{L^\theta(\Omega)}\leq M_1\varepsilon(1+t^{-\frac{3}{2}(\frac{1}{p_0}-\frac{1}{\theta})})e^{-\alpha_1t},\,\forall\theta\in[q_0,\infty];
\\
\!\|\nabla c(\cdot,t)\|_{L^\infty(\Omega)}\leq M_2\varepsilon(1+t^{-\frac12})e^{-\alpha_1 t} \quad \hbox{for all }\, t\in[0,\widetilde{T});\\
\|u(\cdot,t)\|_{L^{q_0}(\Omega)}\leq M_3\varepsilon\left(1+t^{-\frac12+\frac{3}{2q_0}}\right) e^{-\alpha_2 t}
~\hbox{for all}~t\in[0,\widetilde{T});\\
\|\nabla u(\cdot,t)\|_{L^{3}(\Omega)}\leq M_4\varepsilon\left(1+t^{-\frac12}\right) e^{-\alpha_2 t}
~\hbox{for all}~t\in[0,\widetilde{T}).
\end{array}\right.
\!\right\}
\end{align}}
By Lemma 2.5 and (1.6), $T>0$ is well-defined. As in the previous subsection,  we first show  $T=T_{max}$, and then $T_{max}=\infty$. To this end,
we will show that all of the estimates
mentioned in \eqref{3.79} are valid with even smaller coefficients on the right hand side than that in \eqref{3.79}. The derivation of these estimates will mainly rely on $L^p-L^q$ estimates for the Neumann heat semigroup and the corresponding semigroup for Stokes operator, and the fact that the classical solution of \eqref{1.1} on $(0,T)$ can be represented as 
\setlength\belowdisplayskip{4pt}
\begin{align}
(m-\rho)(\cdot,t)\!&=\!e^{t\Delta}(m_0-\rho_0)\!+\!\int_0^te^{(t-s)\Delta}(\nabla\cdot(\rho \mathcal{S}(x,\rho,c)\nabla c)-u\cdot\nabla(m-\rho))(\cdot,s)ds,\label{3.80}
\\
&\rho(\cdot,t)=e^{t\Delta}\rho_0-\int_0^te^{(t-s)\Delta}(\nabla\cdot(\rho \mathcal{S}(x,\rho,c)\nabla c)+u\cdot\nabla\rho+\rho m)(\cdot,s)ds,\label{3.81}
\\
&c(\cdot,t)=e^{t(\Delta-1)}c_0+\int_0^te^{(t-s)(\Delta-1)}(m-u\cdot\nabla c)(\cdot,s)ds,\label{3.82}
\\
u(\cdot,t)=& e^{-tA}u_0+\int_0^te^{-(t-s)A}\mathcal{P}((\rho+m)\nabla\phi-
 (u\cdot\nabla)u )(\cdot,s)ds.\label{3.83}
\end{align}}
The proofs of the following two lemmas  are same as  that of \cite{LPW}, so we omit it here.

\begin{lemma} (Lemma 3.17 in \cite{LPW})\label{lemma3.18}
Under the assumptions of Proposition 3.2,
\begin{align*}
\|(m-\rho)(\cdot,t)-m_\infty\|_{L^\theta(\Omega)}\leq M_5\varepsilon(1+t^{-\frac32(\frac{1}{p_0}-\frac{1}{\theta})})e^{-\alpha_1t}
\end{align*}
for all $t\in(0,T)$ and $\theta\in[q_0,\infty]$ with $M_5=1+C_1+C_1|\Omega|^{\frac1{p_0}-\frac1{q_0}}$.
\end{lemma}

\begin{lemma}(Lemma 3.18 in \cite{LPW})\label{lemma3.19}
Under the assumptions of Proposition 3.2,
{\setlength\abovedisplayskip{3pt}
\setlength\belowdisplayskip{3pt}
\begin{align*}
\|m(\cdot,t)-m_\infty\|_{L^{\theta}(\Omega)}\leq (M_5+M_1)\varepsilon(1+t^{-\frac 32(\frac{1}{p_0}-\frac{1}{\theta})}) e^{-\alpha_1 t}
\quad\hbox{for all}\,\,t\in(0,T), \theta\in[q_0,\infty].
\end{align*}}
\end{lemma}

\begin{lemma}\label{lemma3.20}
Under the assumptions of Proposition 3.2,  we have
\begin{align*}
\|u(\cdot,t)\|_{L^{q_0}(\Omega)}\leq \frac{ M_3}2\varepsilon(1+t^{-\frac12+\frac{3}{2q_0}}) e^{-\alpha_2 t}
\quad ~\hbox{for all }~t\in(0,T).\end{align*}
\end{lemma}

\proof For any given $\alpha_2<\lambda_1'$, we can fix $ \mu\in (\alpha_2, \lambda_1')$.
 By  \eqref{3.83}, Lemma 2.2, Lemma 2.3 and $\mathcal{P}(\nabla \phi)=0$, we obtain that
{\setlength\abovedisplayskip{4pt}
\setlength\belowdisplayskip{4pt}\begin{align}
&\|u(\cdot,t)\|_{L^{q_0}(\Omega)}\nonumber
\\
\leq& C_6t^{-\frac 32(\frac13-\frac{1}{q_0})}e^{-\mu t}\|u_0\|_{L^3(\Omega)}
+\int_0^t\|e^{-(t-s)A}\mathcal{P}((\rho+m)\nabla\phi-(u\cdot \nabla) u)(\cdot,s)\|_{L^{q_0}(\Omega)}ds\label{3.84}
\\
\leq& C_6t^{-\frac 12+\frac{3}{2q_0}}e^{-\mu t}\varepsilon
+C_6C_9\|\nabla\phi\|_{L^\infty(\Omega)}\int_0^te^{-\mu(t-s)}\|(\rho+m-m_\infty)(\cdot,s)\|_{L^{q_0}(\Omega)}ds\nonumber\\
& +C_6C_9 \int_0^t(t-s)^{-\frac 12}e^{-\mu(t-s)}\|(u \cdot\nabla)u(\cdot,s)\|_{L^{\frac 1{\frac 13+\frac 1{q_0}}}(\Omega)}ds. \nonumber
\end{align}}
By Lemma \ref{lemma3.19} and the definition of $T$, we get
{\setlength\abovedisplayskip{4pt}
\setlength\belowdisplayskip{4pt}\begin{align}
\|(\rho+m-m_\infty)(\cdot,s)\|_{L^{q_0}(\Omega)}
=&\|(m-m_\infty)(\cdot,s)\|_{L^{q_0}(\Omega)}+\|\rho(\cdot,s)\|_{L^{q_0}(\Omega)}\label{3.85}
\\
\leq &(2M_5+M_1)\varepsilon(1+s^{-\frac{3}{2}(\frac{1}{p_0}-\frac{1}{q_0})})e^{-\alpha_1s}.\nonumber
\end{align}}
Inserting \eqref{3.85} into \eqref{3.84},  by the definition of $T$
and noting that $\frac{3}{2}(\frac{1}{p_0}-\frac{1}{q_0})<1$, we have 
{\setlength\abovedisplayskip{4pt}
\setlength\belowdisplayskip{4pt}\begin{align*}
&\|u(\cdot,t)\|_{L^{q_0}(\Omega)}
\\
\leq& C_6t^{-\frac 12+\frac{3}{2q_0}}e^{-\mu t}\varepsilon
+C_6C_9(2M_5+M_1)\|\nabla\phi\|_{L^\infty(\Omega)}\varepsilon\int_0^t(1+s^{-\frac 32(\frac{1}{p_0}-\frac{1}{q_0})})e^{-\alpha_1s}e^{-\mu(t-s)}ds
\\
& +C_6C_9 \int_0^t(t-s)^{-\frac 12}e^{-\mu(t-s)} \|\nabla u(\cdot,s)\|_{L^{3}(\Omega)}\|u(\cdot,s)\|_{L^{q_0}(\Omega)}ds \nonumber\\
\leq& C_6 t^{-\frac12+\frac{3}{2q_0}}e^{-\mu t}\varepsilon
+C_6C_9 C_{10}(2M_5+M_1)\|\nabla\phi\|_{L^\infty(\Omega)}\varepsilon(1+t^{\min\{0,1-\frac 32(\frac{1}{p_0}-\frac{1}{q_0})\}})e^{-\alpha_2t}
 \nonumber\\
  & +  3 C_6C_9 M_3 M_4 \varepsilon^2 \int_0^t(t-s)^{-\frac 12}   e^{-\mu(t-s)}
(1+s^{-1+\frac{3}{2q_0}}) e^{-2\alpha_2 s} ds
 \nonumber\\
\leq& C_6 t^{-\frac12+\frac{3}{2q_0}}\varepsilon e^{-\mu t}
+2C_6C_9 C_{10}(2M_5+M_1) \|\nabla\phi\|_{L^\infty(\Omega)}\varepsilon e^{-\alpha_2 t}\nonumber\\
 &+  3 C_6C_9 C_{10}M_3 M_4 (1+t^{-\frac12+\frac{3}{2q_0}})\varepsilon^2e^{-\alpha_2 t}  \nonumber\\
\leq & \frac{ M_3}2\varepsilon(1+t^{-\frac12+\frac{3}{2q_0}}) e^{-\alpha_2 t},
\end{align*}}
where we have used    \eqref{34.1b} and \eqref{34.3}.

\begin{lemma}\label{Lemma 3.14b}
Under the assumptions of Proposition 3.2, we have
$$
\|\nabla u(\cdot,t)\|_{L^{3}(\Omega)}\leq \frac {M_4}2 \varepsilon(1+t^{-\frac12})  e^{-\alpha_2 t}~~\hbox{for all}~~t\in(0,T).
$$
\end{lemma}

\proof According to  \eqref{3.83}, and
applying Lemma 2.2(iii) and Lemma 2.3, we arrive at
{\setlength\abovedisplayskip{4pt}
\setlength\belowdisplayskip{4pt}
\begin{align}
&\| \nabla u(\cdot,t)\|_{L^{3}(\Omega)}\nonumber
\\
\leq& C_7 t^{-\frac 12}e^{-\mu t}\|u_0\|_{L^3(\Omega)}
+\int_0^t\| \nabla e^{-(t-s)A}\mathcal{P}((\rho+m)\nabla\phi- (u \cdot\nabla)u)(\cdot,s)\|_{L^3(\Omega)}ds\nonumber
\\
\leq& C_7 t^{-\frac 12}e^{-\mu t}\varepsilon
+C_7|\Omega|^{\frac1{3}\!-\!\frac1{q_0}}\int_0^t (t-s)^{-\frac{1}{2}}e^{-\mu(t-s)}\|\mathcal{P}((\rho+m-m_{\infty})\nabla\phi)(\cdot,s)\|_{L^{q_0}(\Omega)}ds\nonumber\\
& +C_7\int_0^t  (t-s)^{-\frac 12-\frac{3}{2q_0}}  e^{-\mu(t-s)}\|\mathcal{P}((u \cdot\nabla)u)(\cdot,s)\|_
{L^{\frac {3q_0}{3+q_0}}(\Omega)
}
ds\label{3.58b}
\\
\leq&  C_7 t^{-\frac 12}e^{-\mu t}\varepsilon
+C_7C_9\|\nabla\phi\|_{L^\infty(\Omega)} |\Omega|^{\frac1{3}\!-\!\frac1{q_0}}
\int_0^t  (t-s)^{-\frac{1}{2}} e^{-\mu(t-s)}\|(\rho+m-m_{\infty})(\cdot,s)\|_{L^{q_0}(\Omega)}ds\nonumber\\
& +C_7C_9\int_0^t  (t-s)^{-\frac 12-\frac{3}{2q_0}}  e^{-\mu(t-s)} \|\nabla u(\cdot,s)\|_{L^{3}(\Omega)}\|u(\cdot,s)\|_{L^{q_0}(\Omega)}ds,
 \nonumber
\end{align}
where $\mathcal{P}(m_{\infty}\nabla\phi)=m_{\infty} \mathcal{P}(\nabla\phi)=0$ is used.

From \eqref{3.85}, it  follows that
{\setlength\abovedisplayskip{4pt}
\setlength\belowdisplayskip{4pt}
\begin{align}
& \int_0^t  (t-s)^{-\frac{1}{2}} e^{-\mu(t-s)}\|(\rho+m-m_{\infty})(\cdot,s)\|_{L^{q_0}(\Omega)}ds\label{1.50b}\\
\leq&(2M_5+M_1)\varepsilon
\int_0^t(t-s)^{-\frac{1}{2}}e^{-\mu(t-s)}(1+s^{-\frac{3}{2}(\frac{1}{p_0}-\frac{1}{q_0})})e^{-\alpha_1s}ds. \nonumber
\end{align}
In addition, an  application of the H\"{o}lder inequality and definition of $T$ 
shows  that
\begin{align}
& \int_0^t(t-s)^{-\frac 12-\frac{3}{2q_0}}  e^{-\mu(t-s)}\|u(\cdot,s)\|_{L^{q_0}(\Omega)}
\|\nabla u(\cdot,s)\|_{L^{3}(\Omega)}ds
\nonumber\\
\leq&
3 M_3 M_4 \varepsilon^2  \int_0^t(t-s)^{-\frac 12-\frac{3}{2q_0}}   e^{-\mu(t-s)}
(1+s^{-1+\frac{3}{2q_0}}) e^{-2\alpha_2 s} ds.
\label{1.50c}
\end{align}
Therefore, inserting \eqref{1.50c}, \eqref{1.50b} into \eqref{3.58b} and applying Lemma 2.4, we get
{\setlength\abovedisplayskip{4pt}
\setlength\belowdisplayskip{4pt}
\begin{align*}
&\|\nabla u(\cdot,t)\|_{L^{3}(\Omega)}
\\
\leq& C_7 t^{-\frac12}e^{-\mu t}\varepsilon
+C_7C_9 C_{10} \|\nabla\phi\|_{L^\infty(\Omega)}|\Omega|^{\frac1{3}\!-\!\frac1{q_0}}(2M_5+M_1)\varepsilon(1+t^{\min\{0,\frac12-\frac 3{ 2} (\frac{1}{p_0}-\frac{1}{q_0})\}})e^{-\alpha_2t}\\
& +3 C_7C_9 C_{10} M_3 M_4 \varepsilon^2(1+t^{-\frac12} )e^{-\alpha_2 t}
\\
\leq & \frac {M_4}2 \varepsilon(1+t^{-\frac12}) e^{-\alpha_2 t},
\end{align*}}
where \eqref{34.1c}, \eqref{34.3b} are used.

\begin{lemma}\label{lemma3.21}
Under the assumptions of Proposition  3.2, we have{\setlength\abovedisplayskip{4pt}
\setlength\belowdisplayskip{4pt}
\begin{align*}
\|\nabla c(\cdot,t)\|_{L^{\infty}(\Omega)}\leq \frac{M_2}{2}\varepsilon(1+t^{-\frac12}) e^{-\alpha_1 t}\quad \hbox{for all}\,\, t\!\in\!(0,T).
\end{align*}}
\end{lemma}
\vspace{-1em}
\proof From \eqref{3.82} and the standard regularization properties of the Neumann heat semigroup $(e^{\tau\Delta})_{\tau>0}$ in \cite{Winkler7}, one can conclude that
{\setlength\abovedisplayskip{4pt}
\setlength\belowdisplayskip{4pt}
\begin{align}
\|\nabla c(\cdot,t)\|_{L^\infty(\Omega)}
&\leq e^{-t}\|\nabla e^{t \Delta}c_0\|_{L^\infty(\Omega)}+\int_0^t\|\nabla
e^{(t-s)(\Delta-1)}(m-u\cdot\nabla c)(\cdot,s)\|_{L^\infty(\Omega)}ds\nonumber
\\
&\leq  C_3(1+t^{-\frac12})e^{- t}\|\nabla c_0\|_{L^3(\Omega)}+\int_0^t\|\nabla
e^{(t-s)(\Delta-1)}(m-m_{\infty})(\cdot,s)\|_{L^\infty(\Omega)}ds\nonumber
\\
&\quad+\int_0^t\|\nabla
e^{(t-s)(\Delta-1)}u\cdot\nabla c(\cdot,s)\|_{L^\infty(\Omega)}ds.\label{3.86}
\end{align}}
In the second inequality, we have used $ \nabla e^{(t-s)(\Delta-1)}m_{\infty}=0$.

 From Lemma 2.1(ii), Lemma \ref{lemma3.19} and Lemma 2.4, it follows that
{\setlength\abovedisplayskip{4pt}
\setlength\belowdisplayskip{4pt}
\begin{align}
&\int_0^t\|\nabla
e^{(t-s)(\Delta-1)}(m-m_{\infty})(\cdot,s)\|_{L^\infty(\Omega)}ds\nonumber
\\
\leq&C_2\int_0^t(1+(t-s)^{-\frac12-\frac{3}{2q_0}})e^{-(\lambda_1+1)(t-s)}\|(m-m_{\infty}) (\cdot,s)\|_{L^{q_0}(\Omega)}ds\label{3.87}
\\
\leq&C_2 (M_5+M_1)\varepsilon
\int_0^t(1+(t-s)^{-\frac12-\frac{3}{2q_0}})e^{-(\lambda_1+1)(t-s)}
(1+s^{-\frac 32(\frac{1}{p_0}-\frac{1}{q_0})}) e^{-\alpha_1 s}ds\nonumber
\\
\leq&C_2 C_{10}(M_5+M_1)\varepsilon (1+t^{\min\{0,\frac12-\frac{3}{2p_0}\}})e^{-\min\{\alpha_1,\lambda_1+1\}t}\nonumber
\\
\leq&C_2 C_{10}(M_5+M_1)\varepsilon
(1+t^{-\frac12})e^{-\alpha_1t}.\nonumber
\end{align}}
On the other hand,  by  Lemma 2.1(ii), Lemma 2.4 and the definition of $T$, we obtain
{\setlength\abovedisplayskip{4pt}
\setlength\belowdisplayskip{4pt}
\begin{align}
&\int_0^t\|\nabla
e^{(t-s)(\Delta-1)}u\cdot\nabla c(\cdot,s)\|_{L^\infty(\Omega)}ds\nonumber
\\
\leq&C_2\int_0^t(1+(t-s)^{-\frac12-\frac{3}{2q_0}})e^{-(\lambda_1+1)(t-s)}\|u\cdot\nabla c(\cdot,s)\|_{L^{q_0}(\Omega)}ds\label{3.88}
\\
\leq&C_2\int_0^t(1+(t-s)^{-\frac12-\frac{3}{2q_0}})e^{-(\lambda_1+1)(t-s)}\|u(\cdot,s)\|_{L^{q_0}(\Omega)}\|\nabla c(\cdot,s)
\|_{L^\infty(\Omega)}ds\nonumber
\\
\leq& C_2M_3 M_2\varepsilon^2
\int_0^t(1+(t-s)^{-\frac12-\frac{3}{2q_0}})e^{-(\lambda_1+1)(t-s)}
(1+s^{-\frac12+\frac{3}{2q_0}}) (1+s^{-\frac12}) e^{-(\alpha_1+\alpha_2) s} \nonumber
\\
\leq&3C_2M_3M_2\varepsilon^{2}\int_0^te^{-(\lambda_1+1)(t-s)}
 e^{-(\alpha_1+\alpha_2) s}(1+(t-s)^{-\frac12-\frac{3}{2q_0}})(1+s^{-1+\frac{3}{2q_0}})ds\nonumber
\\
\leq&3C_2M_3M_2C_{10}\varepsilon^{2}(1+t^{-\frac12})e^{-\alpha_1t}. \nonumber
\end{align}}
Hence combining above inequalities  and applying \eqref{34.1} and \eqref{34.2}, we arrive at the conclusion.\vspace{-0.5em}

\begin{lemma}\label{lemma3.22}
Under the assumptions of Proposition 3.2, we have
\begin{align*}
\|\rho(\cdot,t)\|_{L^\theta(\Omega)}\leq \frac{M_1}{2}\varepsilon(1+t^{-\frac{3}{2}(\frac{1}{p_0}-\frac{1}{\theta})})e^{-\alpha_1t}
\quad\hbox{for all}\,\,t\in(0,T),\,\theta\in[q_0,\infty].
\end{align*}
\end{lemma}
\proof  From \eqref{3.81}, we have
\begin{align*}
\rho(\cdot,t)=&e^{t(\Delta-m_\infty)}\rho_0-\int_0^te^{(t-s)(\Delta-m_\infty)}(\nabla\cdot(\rho\mathcal{S}(\cdot,\rho,c)\nabla c)-u\cdot\nabla\rho)(\cdot,s)ds\\
&+
\int_0^te^{(t-s)(\Delta-m_\infty)} \rho (m_\infty-m)(\cdot,s)ds
.
\end{align*}
By Lemma \ref{Lemma 2.1},
the result in Section 2 of \cite{Winkler7} and $\alpha_1<\min\{\lambda_1,m_\infty\}$, we obtain
{\setlength\abovedisplayskip{4pt}
\setlength\belowdisplayskip{4pt}
\begin{align*}
&\|\rho(\cdot,t)\|_{L^\theta(\Omega)}
\\
\leq&
e^{-m_\infty t}(\|e^{t\Delta}(\rho_0-\overline{\rho}_0)\|_{L^\theta(\Omega)}+ \|\overline{\rho}_0 \|_{L^\theta(\Omega)}) +\int_0^t\|e^{(t-s)(\Delta-m_\infty)}\nabla\cdot(\rho \mathcal{S}(\cdot,\rho,c)\nabla c)(\cdot,s)\|_{L^\theta(\Omega)}ds
\\
&+\int_0^t\|e^{(t-s)(\Delta-m_\infty)}(u\cdot\nabla\rho)(\cdot,s)\|_{L^\theta(\Omega)}ds
+\int_0^t\|e^{(t-s)(\Delta-m_\infty)}\rho (m_\infty-m)(\cdot,s)\|_{L^\theta(\Omega)}ds
\\
\leq&
C_1(1+t^{-\frac 32(\frac1{p_0}-\frac1{\theta})})e^{-(\lambda_1+m_\infty)t}\|\rho_0-\overline{\rho}_0\|_{L^{p_0}(\Omega)}
+(\min\{1,|\Omega|\})^{-\frac 1{p_0}}e^{-m_\infty t}\varepsilon \\
& +C_4C_S\int_0^t(1+(t-s)^{-\frac12-\frac{3}{2}(\frac{1}{q_0}-\frac{1}{\theta})})e^{-(\lambda_1+m_\infty)(t-s)}\|\rho\|_{L^{q_0}(\Omega)}\|\nabla c\|_{L^\infty(\Omega)}ds\\
&+\int_0^t\|e^{(t-s)(\Delta-m_\infty) }\nabla\cdot(\rho u)(\cdot,s)\|_{L^\theta(\Omega)}ds
+\int_0^t\|e^{(t-s)(\Delta-m_\infty)}\rho (m_\infty-m)(\cdot,s)\|_{L^\theta(\Omega)}ds\\
\leq&
(2C_1+(\min\{1,|\Omega|\})^{-\frac 1{p_0}}) (1+t^{-\frac 32(\frac1{p_0}-\frac1{\theta})}) \varepsilon e^{-\alpha_1t}
\\
&+C_4C_S\int_0^t(1+(t-s)^{-\frac12-\frac{3}{2}(\frac{1}{q_0}-\frac{1}{\theta})})
e^{-(\lambda_1+m_\infty)(t-s)}\|\rho\|_{L^{q_0}(\Omega)}\|\nabla c\|_{L^\infty(\Omega)}ds
\\
&+C_4\int_0^t(1+(t-s)^{-\frac12-\frac{3}{2}(\frac{1}{q_0}-\frac{1}{\theta})})
e^{-(\lambda_1+m_\infty)(t-s)}\|\rho\|_{L^{\infty}(\Omega)}\|u\|_{L^{q_0}(\Omega)}ds
\\
& +C_1\int_0^t(1+(t-s)^{-\frac{3}{2}(\frac{1}{q_0}-\frac{1}{\theta})})
e^{-m_\infty(t-s)}\|\rho\|_{L^{q_0}(\Omega)}\|m-m_\infty\|_{L^{\infty }(\Omega)}ds.
\end{align*}}
According to the definition of $T$,  Lemma \ref{lemma3.21} and Lemma 2.4, this shows that
\begin{align}
&\int_0^t(1+(t-s)^{-\frac12-\frac{3}{2}(\frac{1}{q_0}-\frac{1}{\theta})})
e^{-(\lambda_1+m_\infty)(t-s)}\|\rho\|_{L^{q_0}(\Omega)}\|\nabla c\|_{L^\infty(\Omega)}ds  \nonumber\\
\leq &3  M_1M_2\varepsilon^2 \int_0^t(1+(t-s)^{-\frac12-\frac{3}{2}(\frac{1}{q_0}-\frac{1}{\theta})})
e^{-\lambda_1(t-s)}e^{-2\alpha_1 s}(1+s^ {-\frac12-\frac{3}{2}(\frac{1}{p_0}-\frac{1}{q_0})})
ds\nonumber
\\
\leq&3C_{10} M_1M_2 \varepsilon^2
(1+t^{\min\{0,-\frac{3}{2}(\frac{1}{p_0}-\frac{1}{\theta})\}})
e^{-\min\{\lambda_1,2\alpha_1\}t}.\nonumber
\end{align}
Similarly, we can also get
{\setlength\abovedisplayskip{4pt}
\setlength\belowdisplayskip{4pt}
\begin{align}
&\int_0^t(1+(t-s)^{-\frac12-\frac{3}{2}(\frac{1}{q_0}-\frac{1}{\theta})})
e^{-(\lambda_1+m_\infty)(t-s)}\|\rho\|_{L^{\infty}(\Omega)}\|u\|_{L^{q_0}(\Omega)}ds\nonumber\\
\leq &
3  M_1M_3\varepsilon^2 \int_0^t(1+(t-s)^{-\frac12-\frac{3}{2}(\frac{1}{q_0}-\frac{1}{\theta})})
e^{-\lambda_1(t-s)}e^{-2\alpha_1 s}(1+s^ {-\frac12-\frac{3}{2}(\frac{1}{p_0}-\frac{1}{q_0})})
ds\nonumber
\\
\leq &3C_{10}M_3 M_1\varepsilon^2(1+t^{\min\{0,-\frac{3}{2}(\frac{1}{p_0}-\frac{1}{\theta})\}})
e^{-\min\{\lambda_1,2\alpha_1\}t}\nonumber
\end{align}
and
\begin{align}
&\int_0^t(1+(t-s)^{-\frac{3}{2}(\frac{1}{q_0}-\frac{1}{\theta})})
e^{-m_\infty(t-s)}\|\rho\|_{L^{q_0}(\Omega)}\|m-m_\infty\|_{L^{\infty }(\Omega)}ds\nonumber\\
\leq &
3  M_1(M_5+M_1)\varepsilon^2 \int_0^t(1+(t-s)^{-\frac12-\frac{3}{2}(\frac{1}{q_0}-\frac{1}{\theta})})
e^{-m_\infty(t-s)}
  e^{-2\alpha_1 s}(1+s^{-\frac 3{p_0}+\frac{3}{2 q_0}})
ds\nonumber
\\
\leq&3C_{10}M_1(M_5+M_1)
\varepsilon^2(1+t^{\min\{0,-\frac{3}{2}(\frac{1}{p_0}-\frac{1}{\theta})\}})
e^{-\min\{m_\infty,2\alpha_1\}t},\nonumber
\end{align}}
where the fact that $q_0\in (3,\frac{3p_0}{2(3-p_0)})$ warrants $-\frac 3{p_0}+\frac{ 3}{2 q_0}>-1$ is used.
Hence the combination of the above inequalities  yields
$\|\rho(\cdot,t)\|_{L^\theta(\Omega)}\leq\frac{M_1}{2}\varepsilon(1+t^{-\frac 32(\frac1{p_0}-\frac1{\theta})})e^{-\alpha_1t}$, thanks to  \eqref{34.6a}, \eqref{34.5a} and \eqref{34.3d}.

\begin{lemma}\label{lemma3.23}
Under the assumptions of Proposition 3.2, we have{\setlength\abovedisplayskip{4pt}
\setlength\belowdisplayskip{4pt}
\begin{align*}
\|(m-\rho)(\cdot,t)-e^{t\Delta}(m_0-\rho_0)\|_{L^\theta(\Omega)}\leq \frac{\varepsilon }{2}(1+t^{-\frac{3}{2}(\frac{1}{p_0}-\frac{1}{\theta})}) e^{-\alpha_1 t}
\,\,\hbox{for }\,\theta\in[q_0,\infty], t\in(0,T).\end{align*}}
\end{lemma}

\proof From \eqref{3.80} and Lemma 2.1(iv), it follows that
{\setlength\abovedisplayskip{4pt}
\setlength\belowdisplayskip{4pt}
\begin{align*}
&\|(m-\rho)(\cdot,t)-e^{t\Delta}(m_0-\rho_0)\|_{L^\theta(\Omega)}
\\
\leq&\int_0^t\|e^{(t-s)\Delta}(\nabla\cdot(\rho \mathcal{S}(\cdot,\rho,c)\nabla c)-u\cdot\nabla(m-\rho))(\cdot,s)\|_{L^\theta(\Omega)}ds
\\
\leq&\int_0^t\|e^{(t-s)\Delta}\nabla\cdot(\rho\mathcal{S}(\cdot,\rho,c)\nabla c)(\cdot,s)\|_{L^\theta(\Omega)}ds
+\int_0^t\|e^{(t-s)\Delta}\nabla\cdot((m-\rho-m_\infty)u)(\cdot,s)\|_{L^\theta(\Omega)}ds
\\
\leq&C_4C_S \int_0^t(1+(t-s)^{-\frac12-\frac{3}{2}(\frac{1}{q_0}-\frac{1}{\theta})})
e^{-\lambda_1(t-s)}\|\rho(\cdot,s)\|_{L^{q_0}(\Omega)}\|\nabla c(\cdot,s)\|_{L^\infty(\Omega)}ds
\\
&+C_4\int_0^t(1+(t-s)^{-\frac12-\frac{3}{2}(\frac{1}{q_0}-\frac{1}{\theta})})
e^{-\lambda_1(t-s)}\|u(m-\rho-m_\infty)(\cdot,s)\|_{L^{q_0}(\Omega)}ds
\\
=:&I_1+I_2.
\end{align*}}
\vspace{-0.5em}
From the definition of $T$ and \eqref{34.4}, we have
\begin{align*}
I_1\leq
&
3C_4C_S M_1 M_2 \varepsilon^2 \int_0^t(1+(t-s)^{-\frac12-\frac{3}{2}(\frac{1}{q_0}-\frac{1}{\theta})})
e^{-\lambda_1(t-s)}
(1+s^{-\frac12-\frac{3}{2}(\frac{1}{p_0}-\frac{1}{q_0})}) e^{-2\alpha_1 s}ds \nonumber
\\
 \leq&3C_4C_SC_{10} M_1M_2 \varepsilon^2
(1+t^{\min\{0,-\frac{3}{2}(\frac{1}{p_0}-\frac{1}{\theta})\}})
e^{-\min\{\lambda_1,2\alpha_1\}t}\nonumber
\\
\leq&
\frac{\varepsilon }{4}(1+t^{-\frac 32(\frac1{p_0}-\frac1{\theta})})e^{-\alpha_1t}.\nonumber
\end{align*}
From  Lemma \ref{lemma3.18}, Lemma \ref{lemma3.20} and \eqref{34.7}, it follows that
\begin{align}
I_2&=C_4\int_0^t(1+(t-s)^{-\frac12-\frac{3}{2}(\frac{1}{q_0}-\frac{1}{\theta})})
e^{-\lambda_1(t-s)}\|m-\rho-m_\infty\|_{L^\infty(\Omega)}\|u\|_{L^{q_0}(\Omega)}ds\nonumber
\\
&\leq C_4M_3
M_5\varepsilon^{2}\int_0^t(1+(t-s)^{-\frac12-\frac{3}{2}(\frac{1}{q_0}-\frac{1}{\theta})})
e^{-\lambda_1(t-s)}\nonumber (1+s^{-\frac{3}{2p_0}})e^{-\alpha_1s}
(1+s^{-\frac12+\frac{3}{2q_0}})e^{-\alpha_2 s}ds\nonumber
\\
&\leq 3C_4M_3
M_5\varepsilon^{2}\int_0^t(1+(t-s)^{-\frac12-\frac{3}{2}(\frac{1}{q_0}-\frac{1}{\theta})})
(1+s^{-\frac12+\frac{3}{2}(\frac{1}{q_0}-\frac{1}{p_0})})\nonumber e^{-\lambda_1(t-s)}e^{-(\alpha_1+\alpha_2) s}ds\nonumber
\\
&\leq 3C_{10}C_4M_3
M_5\varepsilon^{2} e^{-\min\{\lambda_1, \alpha_1+\alpha_2\}t}(1+t^{\min\{0,\frac{3}{2}
(\frac{1}{\theta}-\frac{1}{p_0})\}})\nonumber
\\
&\leq \frac{\varepsilon }{4}(1+t^{-\frac{3}{2}(\frac{1}{p_0}-\frac{1}{\theta})})e^{-\alpha_1 t}.\nonumber
\end{align}
Combining the above inequalities, we  arrive at
$$
\|(\rho-m)(\cdot,t)-e^{t\Delta}(\rho_0-m_0)\|_{L^\theta(\Omega)}\leq \frac{\varepsilon}{2}(1+t^{-\frac{3}{2}(\frac{1}{p_0}-\frac{1}{\theta}
)}) e^{-\alpha_1 t}
$$
 and thus complete the proof of this lemma.

By the above lemmas, one can see that $T=T_{max}$, and the further estimates of solutions are needed to ensure  $T_{max}=\infty$. \vspace{-0.5em}
\begin{lemma}\label{lemma3.24}
Under the assumptions of Proposition 3.2, for all $\beta\in (\frac 34,\min\{\frac 54-\frac{3}{2q_0},1\})$ there exists $M_6>0$ such that
{\setlength\abovedisplayskip{4pt}
\setlength\belowdisplayskip{4pt}\begin{align*}
\|A^\beta u(\cdot,t)\|_{L^2(\Omega)}\leq \varepsilon M_6e^{-\alpha_2 t}
\quad\hbox{for}\,\, t\in(t_0,T_{max})\,\hbox{ with}\, \,t_0=\min\{\frac{T_{max}}6,1\}.\end{align*}}
\end{lemma}

\proof The proof is similar to that of \eqref{3.66}, and thus is omitted here.

\begin{lemma}\label{lemma3.25}
Under the assumptions of Proposition 3.2, there exists $M_7>0$ such that
$
\|c(\cdot,t)-m_\infty\|_{L^{\infty}(\Omega)}\leq M_7 e^{-\alpha_2 t}
$
for all $(t_0,T_{max})$ with $ t_0=\min\{\frac{T_{max}}6,1\}$.
\end{lemma}

\proof We refer the readers to the proof of Lemma 3.24 in \cite{LPW}.

At this position, we can show the proof of   Theorem 1.2 in the case $\mathcal{S}=0$ on $\partial\Omega$.

{\bf Proof of Proposition 3.2.}~
We first show that the solution is global, i.e. $T_{max}=\infty$. To this end,
 according to the extensibility criterion in Lemma \ref{lemma2.5},
it suffices to show that there exists $C>0$ such that for all $ t_0<t<T_{max}$ {\setlength\abovedisplayskip{4pt}
\setlength\belowdisplayskip{4pt}
  $$ \|\rho(\cdot,t)\|_{L^\infty(\Omega)}+\|m(\cdot,t)\|_{L^\infty(\Omega)}
+\|c(\cdot,t)\|_{W^{1,\infty}(\Omega)}+
\|A^{\beta}u(\cdot,t)\|_{L^2(\Omega)}<C.$$}
From  Lemma \ref{lemma3.19}, Lemma \ref{lemma3.22} and Lemma \ref{lemma3.24},  there exists  $K_i>0$, $i=1,2,3,4$, such that
{\setlength\abovedisplayskip{4pt}
\setlength\belowdisplayskip{4pt}
\begin{align*}
\|m(\cdot,t)-m_\infty\|_{L^{\infty}(\Omega)}&\leq K_1 e^{-\alpha_1 t},~~
\|\rho(\cdot,t)\|_{L^\infty(\Omega)}\leq K_2e^{-\alpha_1 t},
\end{align*}}
\begin{align*}
\|\nabla c(\cdot,t)\|_{L^{\infty}(\Omega)}\leq  K_3 e^{-\alpha_1 t}, ~~\|A^{\beta}u(\cdot,t)\|_{L^2(\Omega)}
\leq K_4 e^{-\alpha_2 t}
\end{align*}
for $t\in(t_0,T_{max})$.
Furthermore,
 Lemma \ref{lemma3.25}  implies that
$\|c(\cdot,t)-m_\infty\|_{W^{1,\infty}(\Omega)}\leq K_3'e^{-\alpha_2t}$
with some $K_3'>0$ for all $t\in(t_0,T_{max})$.
Since $D(A^\beta)\hookrightarrow L^\infty(\Omega)$ with $\beta\in(\frac{3}{4},1)$, it  follows from Lemma \ref{lemma3.24} that
$
\|u(\cdot,t)\|_{L^\infty(\Omega)}\leq K_4 e^{-\alpha_2 t}
$
for some $K_4>0$ for all $t\in(t_0,T_{max})$. This completes the proof of Proposition 3.2.

\section{Proof of main results for general $\mathcal{S}$}\vspace{-1em}

In this section, we give the proof of our results for the general matrix-valued $\mathcal{S}$ by a rather standard argument,
which is accomplished by an approximation procedure (see \cite{Cao1} for example). In order to make the previous results applicable,
we introduce a family of smooth functions
$\rho_\eta\in C_0^\infty(\Omega)$ and $0\leq\rho_\eta(x)\leq1$ for $\eta\in(0,1),$ $\lim_{\eta\to0}\rho_\eta(x)=1$
and
let
$\mathcal{S}_\eta(x,\rho,c)=\rho_\eta(x)\mathcal{S}(x,\rho,c).$
Using this definition, we  regularize  \eqref{1.4} as follows
 \begin{equation}\label{4.1}
\left\{
\begin{array}{ll}
(\rho_\eta)_t+u_\eta\cdot\nabla\rho_\eta=\Delta\rho_\eta-\nabla\cdot(\rho_\eta \mathcal{S}_\eta(x,\rho_\eta,c_\eta)\nabla c_\eta)-\rho_\eta m_\eta,
\\
(m_\eta)_t+u_\eta\cdot\nabla m_\eta=\Delta m_\eta-\rho_\eta m_\eta,
\\
(c_\eta)_t+u_\eta\cdot\nabla c_\eta=\Delta c_\eta-c_\eta+m_\eta,
\\
(u_\eta)_t+(u_\eta\cdot\nabla) u_\eta=\Delta u_\eta-\nabla P_\eta+(\rho_\eta+m_\eta)\nabla\phi,\quad\nabla\cdot u_\eta=0,\\
\displaystyle\frac{\partial \rho_\eta}{\partial\nu}=\frac{\partial m_\eta}{\partial\nu}
=\frac{\partial c_\eta}{\partial\nu}=0,~ u_\eta=0
\end{array}
\right.
\end{equation}
with the initial data
{\setlength\abovedisplayskip{4pt}
\setlength\belowdisplayskip{4pt}
\begin{align}
\rho_\eta(x,0)=\rho_0(x),~m_\eta(x,0)=m_0(x),~c(x,0)=c_0(x),~\hbox{and}~u_\eta(x,0)=u_0(x),\quad x\in\Omega.
\label{4.2}
\end{align}}
It is observed that $\mathcal{S}_\eta$ satisfies the additional condition $\mathcal{S}=0$ on $\partial\Omega$.
Therefore based on the discussion in Section 3, under the assumptions of Theorem 1.1 and Theorem 1.2,  problem (4.1)--(4.2) admits a global classical solution $(\rho_\eta,m_\eta,c_\eta,u_\eta, P_\eta)$ that satisfies 
{\setlength\abovedisplayskip{4pt}
\setlength\belowdisplayskip{4pt}
\begin{align*}
\|m_\eta(\cdot,t)-m_\infty\|_{L^\infty(\Omega)}\leq K_1e^{-\alpha_1 t},\quad
\|\rho_\eta(\cdot,t)-\rho_\infty\|_{L^\infty(\Omega)}\leq K_2e^{-\alpha_1 t},
\\
\|c_\eta(\cdot,t)-m_\infty\|_{W^{1,\infty}(\Omega)}\leq K_3e^{-\alpha_2t}, \quad
\|A^\beta u_\eta(\cdot,t)\|_{L^2(\Omega)}\leq K_4 e^{-\alpha_2 t}
\end{align*}}
for some constants $K_i$, $i=1,2,3,4$, and  all $t\geq 0$. Applying a standard procedure such as in Lemma 5.2 and Lemma 5.6
of \cite{Cao1}, one can  obtain a subsequence of $\{\eta_j\}_{j\in \mathbb{N}}$ with $\eta_j\to 0$ as $j\to \infty$ such that
$
\rho_{\eta_j}\rightarrow \rho, ~m_{\eta_j}\rightarrow m, ~c_{\eta_j}\rightarrow c, u_{\eta_j}\rightarrow u \quad
\hbox{in}~ C_{loc}^{\nu,\frac{\nu}2}(\overline\Omega\times (0,\infty))
$
 as $j\rightarrow \infty$ for some $\nu\in (0,1)$.
Moreover, by the arguments as in Lemma 5.7, Lemma 5.8 of \cite{Cao1}, one  can also  show that $(\rho,m,c,u, P)$ is a classical solution of \eqref{1.4}
with the decay properties asserted in Theorem 1.1 and Theorem 1.2, respectively. The proof of our main results is thus complete.
\vspace{-1em}


\end{document}